\documentclass[11pt]{article}
\usepackage{geometry}                
\usepackage{graphicx}
\usepackage{amssymb}
\usepackage{epstopdf}
\usepackage{amsmath,amsthm,amscd,amssymb}
\usepackage{latexsym}
\numberwithin{equation}{section}

\theoremstyle{plain}
\newtheorem{theorem}{Theorem}[section]

\newtheorem{conjecture}[theorem]{Conjecture}

\theoremstyle{definition}
\newtheorem{definition}[theorem]{Definition}

\theoremstyle{remark}
\newtheorem{remark}[theorem]{Remark}

\newtheorem{case[theorem]}{Case}

\def\R{\mathbb R}

\def\diam{\text{diam}}

\title{Theory of dimension for large discrete sets and applications}
\author{Alex Iosevich, Misha Rudnev and Ignacio Uriarte-Tuero}

\begin{document}

\maketitle

\abstract{We define two notions of discrete dimension based on the Minkowski and Hausdorff dimensions in the continuous setting. After proving some basic results illustrating these definitions, we apply this machinery to the study of connections between the Erd\H os and Falconer distance problems in geometric combinatorics and geometric measure theory, respectively.}

\tableofcontents

\vskip.125in

\section{Introduction and statement of main results}\label{IntroductionStatementResults}

\vskip.125in

In this paper we study the notion of dimension for a large finite subset $A$ of ${\mathbb R}^d$, $d \ge 2$, of cardinality $N$, discrete and 1-separated in the sense that $|a-a'| \ge 1$ for all $a \not=a' \in A$. The notion of dimension is well developed in the ``continuous" setting.

\begin{definition} \label{minkowski}
Given $E \subset {[0,1]}^d$ and $\delta>0$, let $N_{\delta}$ denote the smallest possible number of balls of radius $\delta$ needed to cover $E$. If
$$ -\limsup_{\delta \to 0} \frac{\log(N_{\delta})}{\log(\delta)}=
-\liminf_{\delta \to 0} \frac{\log(N_{\delta})}{\log(\delta)},$$ we call the resulting number the Minkowski dimension of $E$, denoted by $\dim_{{\cal M}}(E)$.
\end{definition}

\begin{definition} \label{hausdorff}
Let $E \subset {[0,1]}^d$. Define the Hausdorff dimension of $E$, denoted by $\dim_{{\cal H}}(E)$ to be
$$ \inf \left\{s \ge 0: \mathcal {H}^s_{\infty}(E)=0 \right\},$$
where $$\mathcal {H}^s_{\infty}(E) = \inf \left\{ \sum_i r_i^s : E \subset \bigcup_{i} B(x_i,r_i) \right\},$$ i.e. the infimum is taken over all the possible coverings of $E$ by balls $B(x_i,r_i)$ of centers $x_i$ and radius $r_i$. \end{definition}

One can check that the Hausdorff dimension always exists, while the Minkowski dimension may not, and that $\dim_{{\cal M}}(E) \ge \dim_{{\cal H}}(E)$. Indeed, if $E$ is any countable set, one can easily check that $\dim_{{\cal H}}(E)=0$, whereas $\dim_{{\cal M}}(E)$ may well be positive. For example, if for $a>1$ one defines
$$ E=\left\{n^{-\frac{1}{a}}: n=1,2 \dots \right\}, $$ then one can check by a direct calculation that
$\dim_{{\cal M}}(E)=\frac{a}{1+a}$. For a detailed description of the beautiful mathematics related to the Minkowski and Hausdorff dimension, see, for example, treatises by Mattila \cite{mattila} and Falconer \cite{falconerbook}.

We will later (see section \ref{BasicDefinitionsTheorems}) define a notion of Minkowski and Hausdorff dimension for discrete sets of large cardinality $N$.
More precisely, we will state results about families of sets $A_N \subset {\mathbb R}^d$, so that the cardinality $\# A_N = N$, where $N \to \infty$, and the corresponding Minkowski and Hausdorff dimensions will be denoted as $\dim_{\mathcal{M}}(A_N)$ and $\dim_{\mathcal{H}}(A_N)$ (there should be no confusion since the context should make it clear when we refer to the continuous or the discrete version of these dimensions.) We will also develop in section \ref{BasicDefinitionsTheorems} some basic facts about such a theory of dimension for large discrete sets.

A main application of such machinery is to the study of connections between the Erd\H os and Falconer distance problems in geometric combinatorics and geometric measure theory, respectively. Let us remind the reader what these conjectures say.

\begin{conjecture}\label{ErdosConjecture}\textbf{[Erd\H os distance conjecture]}
Let $A \subset {\mathbb R}^d$, $d \ge 2$, and $\# A = N$, then
$$ \# \Delta(A) \gtrapprox {(\# A)}^{\gamma}, $$
where $\gamma$ can be taken to be $\frac{2}{d}$.
\end{conjecture}

Above, $\Delta(A)$ denotes the distance set
$$ \Delta(A)=\{|a-a'|: a,a' \in A\},$$ with
$$ {|x|}^2=x_1^2+x_2^2+\dots+x_d^2,$$
and $X \lessapprox Y$  ($X \gtrapprox Y$) with the controlling parameter $N$ if for every $\epsilon>0$ there exists $C_{\epsilon}>0$ such that $X \leq C_{\epsilon}N^{\epsilon}Y$ ($X \geq C_{\epsilon}N^{-\epsilon}Y$). If the above notations also allow $\epsilon=0$, we write $X \lesssim Y$  ($X \gtrsim Y$) instead, as well as $X\approx Y$ whenever $X\lesssim Y$ and $X\gtrsim Y$.

Taking $A={[0,N^{\frac{1}{d}}]}^d \cap {\mathbb Z}^d$ shows that one cannot in general do better. In the continuous setting, the analogous conjecture is

\begin{conjecture}\label{FalconerConjecture}\textbf{[Falconer distance conjecture]}
Let $E \subset {[0,1]}^d$ be such that its Hausdorff dimension satisfies $\dim_{\mathcal{H}} (E) > s_0$. Then the Lebesgue measure of $\Delta(E)$ is positive (i.e. $\mathcal{L}^1 (\Delta (E)) >0$.) More precisely, Falconer's conjecture is that $s_0 = \frac{d}{2}$.
\end{conjecture}


Once again taking $E$ to be a set built on an appropriately scaled version of the integer lattice shows that it is possible for $\Delta(E)$ to have Lebesgue measure $0$ if the Hausdorff dimension of $E$ is any number less than [or equal if $d=2$] $\frac{d}{2}$.

See e.g. \cite{pachagarwall} for a thorough description of Erd\H os' conjecture and related problems. Both conjectures have attracted substantial and deep work, and both are far from being proved for any $d$. The best results to date for Erd\H os' conjecture are due to Katz and Tardos \cite{katztardos} in $\R^2$ ($\gamma \approx  0.86$ instead of $1$) and Solymosi and Vu in $d \geq 3$ \cite{solymosivu} ($\gamma$ close to $\frac{2}{d}-\frac{1}{d^2}$.) An earlier result by Solymosi and T{\'o}th \cite{solymositoth}, obtained $\gamma = \frac{6}{7}$ in $\R^2$.

With respect to Falconer's conjecture, after results by Falconer \cite{falconerdistancesets}, Mattila \cite{mattiladistancesets}, and Bourgain \cite{bourgaindistancesets}; Wolff \cite{wolffdistancesets} obtained the best result to date in $\R^2$, namely $s_0 = \frac{4}{3}$, and Erdo{\~g}an \cite{erdogan}, in $d \geq 3$, proved $s_0 = \frac{d}{2}+\frac{1}{3}$.

Work of Katz and Tao, e.g. \cite{katztaoarithmeticprojectionskakeya}, suggests a strong connection between Falconer's conjecture and the Kakeya conjecture (that if $E \subset \R^d$ contains a unit line segment in every direction, then $\dim (E) =d$.)

A full rigorous connection between Erd\H os' and Falconer's conjectures has so far not been established, to our knowledge. The connection between the putative sharpness examples in the Erd\H os and Falconer distance problems led the first named author and I. Laba \cite{iosevichlabakdistance} to prove that the Erd\H os distance conjecture in the special case of Delone sets (which appear naturally in crystallography and in the context of spectral sets in Fourier analysis) is, in fact, implied by the Falconer conjecture (see section \ref{ApplicationsToErdosFalconer}.) Delone sets are roughly speaking statistical perturbations of the integer lattice $\mathbb{Z}^d \subset \mathbb{R}^d$ (see section \ref{ApplicationsToErdosFalconer} for the precise definition.)

One of the main threads of this paper is to further the understanding of such a ``Falconer-to-Erd\H os dictionary", i.e. assuming results of Falconer type, deduce results of Erd\H os type. In this direction, we get the following Theorems \ref{HausdorffAlphaAdaptableSatisfyErdosIfAlphaAccordingToFalconerVersion2WhenSetsAreNestedVersionIntroduction} and
\ref{HausdorffAlphaAdaptableSatisfyErdosIfAlphaAccordingToFalconerVersion2WhenFalconerIsStrengthenedVersionIntroduction}.
Both results essentially state that if the Falconer conjecture holds for dimensions $\alpha > s_0$, then the Erd\H os distance conjecture holds for exponent $\gamma = \frac{1}{s_0}$. However the first result (Theorem \ref{HausdorffAlphaAdaptableSatisfyErdosIfAlphaAccordingToFalconerVersion2WhenSetsAreNestedVersionIntroduction}) assumes the Falconer conjecture as stated, but then has to assume the nesting of the sets $ C_{N,\alpha}={({[diam(A_N)]}^{-1}A_N)}_{N^{-\frac{1}{\alpha}}}$ (which are a fattening by $N^{-\frac{1}{\alpha}}$ of the sets ${[diam(A_N)]}^{-1}A_N $, where given a real number $t>0$, $ tA=\{ta: a \in A\}.$) In turn, the second result (Theorem \ref{HausdorffAlphaAdaptableSatisfyErdosIfAlphaAccordingToFalconerVersion2WhenFalconerIsStrengthenedVersionIntroduction}) does not assume nesting, but has to assume a slightly stronger version of the Falconer conjecture, namely that not only the distance set $\Delta(E)$ has positive length, but that there is a quantitative control of the length $\mathcal{L}^1 (\Delta(E)) \geq C=C(\alpha, C_0)>0$. It should be noted however, that all known recent proofs of results pertaining to the Falconer conjecture actually yield such a quantitative control of the length $\mathcal{L}^1 (\Delta(E))$.

Section \ref{BasicDefinitionsTheorems} contains  precise definitions involved in the formulation of the following theorems. We have nevertheless chosen to give the theorems' formulations here at the outset, assuming that the reader is familiar with the foundations of the continuous dimension theory in terms of energy and capacity (the $\alpha$-energy integral of a measure $\mu$ further being denoted as  $I_\alpha(\mu)$), see Definition \ref{energy}. The discrete version developed further rests on the notion of Hausdorff $\alpha$-adaptability for the families $A_N$ of discrete sets, representing the direct equivalent of the energy condition for the continuous Hausdorff dimension $\alpha$, see Definition \ref{DefinitionHausdorffAlphaAdaptability}. The notion of discrete Hausdorff dimension $\dim_{\mathcal{H}} (A_N) \geq \alpha$, however, only requires $\alpha$-adaptability of ``sufficiently large" subsets $B_N$ of the sets $A_N$. See Definition \ref{DiscreteHausdorffDimension}. In this sense, determining the discrete Hausdorff dimension of a given family of sets is a major problem, alike to that of determining the classical Hausdorff dimension of continuous sets.

Our main results are as follows.

\begin{theorem} \label{HausdorffAlphaAdaptableSatisfyErdosIfAlphaAccordingToFalconerVersion2WhenSetsAreNestedVersionIntroduction}
Suppose that the Falconer distance conjecture holds to the extent that if the Hausdorff dimension of $E \subset {[0,1]}^d$ is greater than $s_0$ ($s_0 \ge \frac{d}{2}$), then the Lebesgue measure of $\Delta(E)$ is positive.
Let $A_N$ be a family of sets with $\# (A_N) = N$ which is Hausdorff $\alpha_0$-adaptable, for some $\alpha_0 > s_0$. Assume also that for any $s_0 < \alpha < \alpha_0$, the family $ C_{N,\alpha}={({[diam(A_N)]}^{-1}A_N)}_{N^{-\frac{1}{\alpha}}}$ is a nested family of sets, i.e. $C_{N+1,\alpha} \subseteq C_{N,\alpha}$.
Then
$$ \# \Delta(A_N) \gtrapprox N^{\frac{1}{s_0}}.$$
\end{theorem}

As far as the nesting requirement in the above theorem is concerned, $N$ can certainly be only a subsequence of integers, in which case $N$ itself in the estimates should be substituted by $\#(A_N)$. We also get another version of Theorem \ref{HausdorffAlphaAdaptableSatisfyErdosIfAlphaAccordingToFalconerVersion2WhenSetsAreNestedVersionIntroduction} under some conditions that are more restrictive than the condition that $\dim_{\mathcal{H}}(A_N) \geq \alpha_0$ (because of a nesting requirement for the ``large subsets" of $A_N$.) Our main Theorem is

\begin{theorem} \label{HausdorffAlphaAdaptableSatisfyErdosIfAlphaAccordingToFalconerVersion2WhenFalconerIsStrengthenedVersionIntroduction}

Suppose that the Falconer distance conjecture holds to the extent that if a Borel probability measure $\mu$ supported on $E \subset {[0,1]}^d$ satisfies that $I_\alpha (\mu) \leq C_0 < \infty$, for some $\alpha > s_0 \ge \frac{d}{2}$ (see Theorem \ref{EnergyCharacterizationOfHausdorffDimension} below), then $\mathcal{L}^1 (\Delta(E)) \geq C=C(\alpha, C_0)>0$.

Let $A_N \subset \mathbb{R}^d$ be a family of sets with $\# (A_N) = N$ with $\dim_{\mathcal{H}}(A_N) = \alpha_0 > s_0$. Then
$$ \# \Delta(A_N) \gtrapprox N^{\frac{1}{s_0}}.$$

In particular, if the Falconer conjecture is true (with the above quantitative control $\mathcal{L}^1 (\Delta(E)) \geq C=C(\alpha, C_0)>0$), then the Erd\H os conjecture is true for any family of sets $A_N \subset \mathbb{R}^d$ with (discrete) Hausdorff dimension $\dim_{\mathcal{H}}(A_N) > \frac{d}{2}$.



\end{theorem}

To better understand the scope of these results, notice first that our Theorem includes the aforementioned result by the first named author and I. Laba \cite{iosevichlabakdistance} (quoted below as Theorem \ref{iosevichlaba}), since we get that

\begin{theorem} \label{DeloneSetsAreDAdaptableVersionIntroduction}
Delone sets in $\mathbb{R}^d$ have discrete Hausdorff dimension $d$.
\end{theorem}

Actually, the class of sets with discrete Hausdorff dimension $\geq \alpha$ is a pretty large class of sets, since, given any set $E \subset \mathbb{R}^d$, of (continuous, i.e. the usual) Hausdorff dimension $\alpha_0$, then for any $\alpha < \alpha_0$, we can build a sequence of sets $A_N$ which is Hausdorff $\alpha$-adaptable, and hence has discrete Hausdorff dimension $\geq \alpha$ (and which, in a sense to be made precise later, ``converges" to (a subset of) $E$.) This is the content of

\begin{theorem}\label{FeketePointsApproximationVersionIntroduction}
Let $E \subset [0,1]^d$ be a compact set so that there exists a Borel probability measure $\mu$ supported on $E$ with $I_\alpha (\mu) < \infty$ (see Theorem \ref{EnergyCharacterizationOfHausdorffDimension}), for $0< \alpha <d$. Then there exists a family of Hausdorff $\alpha$-adaptable sets $A_{N_j} \subset [0,1]^d $, and hence with discrete Hausdorff dimension $\geq \alpha$, with $\# (A_{N_j} ) = N_j \to \infty$, so that, with the notation of \eqref{ProbabilityMeasureMuForAN}, $\mu_{A_{N_j}} \rightharpoonup \mu_0$ (weak-$\ast$ convergence) with $\mu_0$ a Borel probability measure supported on $K_0$ satisfying $I_\alpha (\mu_0) < \infty$, and $A_{N_j} \to \widetilde{K_0}$ in the Hausdorff metric, with $K_0 \subseteq \widetilde{K_0}  \subseteq E$.
\end{theorem}

Regrettably, there is also a class of discrete sets to which the machinery developed does not apply in order to yield results of Erd\H os type. More precisely, the machinery does not apply to families of discrete sets $A_N$ with discrete Hausdorff dimension $\alpha < \frac{d}{2}$ in $\mathbb{R}^d$, since Falconer's conjecture says nothing about such dimensions. However, it should be noted that the techniques from geometric combinatorics allow us to pass from the family of sets $A_N$ to a family of subsets $B_N \subseteq A_N$, provided that the sets $B_N$ are ``sufficiently large" (see section \ref{BasicDefinitionsTheorems} for the precise definitions.) This is why in the definition of discrete Hausdorff dimension we allow also for families of subsets to be taken into account. This allowance for families of subsets sometimes gives rise to surprises. Namely, some families of sets $A_N$ that are not Hausdorff $\alpha$-adaptable for any $\alpha > \frac{d}{2}$ in $\mathbb{R}^d$ (i.e. they would not have discrete Hausdorff dimension $> \frac{d}{2}$ if the families of subsets were not allowed towards computing the discrete Hausdorff dimension), actually ``hide" inside them small copies of ``full dimension" sets, and then the machinery applies to yield for those sets $A_N$ the same kind of Erd\H os type results one would get if the whole sets $A_N$ were ``full dimension" sets (i.e. dimension $d$ in $\mathbb{R}^d$.) Consequently, the class of discrete sets to which the machinery developed does not apply is smaller than what one might think at first sight. That is the content of the example stated below as Theorem \ref{ExampleMinkowskiDimensionDifferentToHausdorffDimension}.


However, we also found families of sets $A_N$ with small Hausdorff dimension (i.e. neither them nor ``hidden" families of sufficiently large subsets $B_N$ are Hausdorff $\alpha$-adaptable for $\alpha$ large). That is the content of Theorems \ref{CantorSetSmallDimension} and \ref{CantorTypeSetDimensionZero} below.
This basically shows that the direct connection between the Erd\H os and Falconer conjectures breaks down whenever the discrete dimension is smaller than $\frac{d}{2}.$

In our opinion, one of the merits of this paper is not so much the techniques we used, which are known in the areas of geometric combinatorics, potential theory and geometric measure theory, but how these techniques and these areas are related in ways not known before to yield the results and ideas we present.

The paper is structured as follows. In section \ref{BasicDefinitionsTheorems} we give the precise basic definitions of the theory of dimension for discrete sets and prove some of the basic Theorems for the understanding of this theory. In section \ref{ApplicationsToErdosFalconer} we give the applications of this machinery to problems of Erd\H os and Falconer type. In section \ref{Examples} we give examples related to the theory.



\vskip.125in

\section{Basic Definitions and Theorems}\label{BasicDefinitionsTheorems}

\vskip.125in

In view of the classical definitions of Minkowski and Hausdorff dimension, how should one define a notion of dimension for discrete sets? The first reasonable step is to control the largest scale by replacing a discrete, one-separated set $A$ of cardinality $N$ by
${[diam(A)]}^{-1}A$, where $diam(A)$ is the diameter of $A$ and given a real number $t>0$,
$$ tA=\{ta: a \in A\}.$$

In order to make a connection with the continuous setting, let us now replace ${[diam(A)]}^{-1}A$ by
${({[diam(A)]}^{-1}A)}_{\delta}$, where given a set $S$, $S_{\delta} = \{ x\in\R^d : d(x,S)\leq \delta \}$ denotes the $\delta$-neighborhood of $S$. If we do not want these $\delta$-balls to interact, we should impose a condition that
$$\delta < \frac{1}{2} \frac{1}{diam(A)}.$$

A discrete variant of the Minkowski dimension now becomes apparent. If after the above procedure  $\delta \approx  \frac{1}{diam(A)}$ happens to be $\delta \approx N^{-\frac{1}{\alpha}}$, $\alpha>0$, $A$ should be a set of Minkowski dimension $\alpha$ (since it is covered by $N$ disjoint balls of radius $\delta$ and $N \delta^\alpha \approx 1$ .) At this point the reader may rightfully point out that
$$ {({[diam(A)]}^{-1}A)}_{N^{-\frac{1}{\alpha}}}$$ has positive Lebesgue measure. However, its measure goes to $0$ as $N$ tends to infinity. The set is, however, uniformly $\alpha$ dimensional in the following sense.

\begin{definition} \label{uniformminkowski} Let $E_N \subset {[0,1]}^d$ be a family of sets dependent on a parameter $N$. Suppose that there exist finite positive constants $C,c$, independent of $N$, such that
$$ c \leq \liminf_{\delta \to 0} \frac{|{(E_N)}_{\delta}|}{\delta^{d-\alpha}} \leq
\limsup_{\delta \to 0} \frac{|{(E_N)}_{\delta}|}{\delta^{d-\alpha}} \leq C,$$ where given a set $S$, $|S|$ denotes its Lebesgue measure. Then we say that the family $E_N$ is uniformly Minkowski $\alpha$-dimensional.
\end{definition}

For the analogy with the continuous case, see e.g. \cite{mattila} p.79.


\begin{theorem} \label{uniformminkowskiest} Let the parameter $N$ run over a subsequence of the natural numbers. Let $A_N \subset {\mathbb R}^d$ be a family of $1$-separated finite sets so that the cardinality of $A_N = \#\{A_N \} =N $. Assume that $${({[diam(A_N)]}^{-1}A_N)}_{\frac{1}{4diam(A_N)}} \subset {[0,1]}^d.$$
Suppose that
\begin{equation} \label{diameterlim} diam(A_N) \lesssim N^{\frac{1}{\alpha}}, \ \ \text{ i.e. that } \ \  diam(A_N) \leq  C N^{\frac{1}{\alpha}}, \end{equation}
with $C$ independent of $N$.

Then ${({[diam(A_N)]}^{-1}A_N)}_{\frac{1}{4C} N^{-\frac{1}{\alpha}}}$ is uniformly Minkowski $\alpha$-dimensional.
\end{theorem}


\begin{proof}

For $\delta = \frac{1}{4C} N^{- \frac{1}{\alpha}}$, we have that

\begin{equation}
\frac{|{({[diam(A_N)]}^{-1}A_N)}_{\delta}|}{\delta^{d-\alpha}} \approx \frac{ N \ (N^{- \frac{1}{\alpha}})^d }{N^{-\frac{d}{\alpha} +1 }} =1
\end{equation}

\end{proof}

This will lead us to a definition of discrete Minkowski dimension. Before that, let us give the following
\begin{definition} \label{discreteminkowski}  Let $A_N \subset {\mathbb R}^d$ be a family of $1$-separated sets, so that the cardinality of $A_N = \#\{A_N \} =N $. Assume that $${({[diam(A_N)]}^{-1}A_N)}_{\frac{1}{4diam(A_N)}} \subset {[0,1]}^d.$$
We say that $A_N$ is adaptable to the discrete Minkowski dimension $\alpha>0$ (or Minkowski $\alpha$-adaptable) if
(\ref{diameterlim}) holds. \end{definition}

The essence of the definition, in view of Theorem \ref{uniformminkowskiest} is that as long as the diameters of our discrete sets are not too large, we can turn them into a set of Minkowski dimension $\alpha>0$ in a canonical way. Since for the discrete Hausdorff dimension (to be defined later) we will allow families of subsets, in order that certain properties remain consistent with the continuous Minkowski and Hausdorff dimensions, we will also allow for subsets here.
\begin{definition} \label{DiscreteMinkowskiDimension}
We define the discrete Minkowski dimension of a family of $1$-separated sets $A_N \subset {\mathbb R}^d$ with $\#\{A_N \} =N $ to be
\begin{eqnarray}
\dim_{\mathcal{M}} (A_N) = \sup \{ \beta >0 & : & \text{ for every } \varepsilon>0, \text{ there exists a family of sets } B_N \subseteq A_N \nonumber \\
& &  \text{ and a constant } C_\varepsilon >0, \text{ so that } \#(B_N) \geq \frac{C_\varepsilon}{N^\varepsilon} \#(A_N), \nonumber \\
& &  \text{ and so that } B_N  \text{ is Minkowski $\beta$-adaptable. } \} \nonumber
\end{eqnarray}
The constant $C_\varepsilon$ depends on $\varepsilon$ and on the sequence $\{B_N \}$, but not on $N$. If there are no such $\beta >0$, the Minkowski dimension of $A_N$ is zero.

\end{definition}

The situation turns out to be far more fascinating with the Hausdorff dimension. We start out by reminding the reader of a connection between the Hausdorff dimension and upper bounds on energy integrals.

\begin{definition} \label{energy} Given a Borel probability measure $\mu$ supported on $E \subset {[0,1]}^d$, the $\alpha$-energy of $\mu$ is given by
$$ I_{\alpha}(\mu)=\int \int {|x-y|}^{-\alpha} d\mu(x)d\mu(y).$$
\end{definition}

A classical result in geometric measure theory connecting energies and dimension is the following (see e.g. \cite{mattila} pp.109-114.)
\begin{theorem}\label{EnergyCharacterizationOfHausdorffDimension}
Let $\alpha$ be the Hausdorff dimension of $E \subset {[0,1]}^d$ and let $\mu$ be a Borel probability measure supported on $E$. Then
$$ \alpha=\sup \left\{s>0: \exists \mu \text{ with } I_s(\mu)<\infty \right\}.$$
\end{theorem}

This leads us to explore the energy integral associated with the Lebesgue measure on
${({[diam(A)]}^{-1}A)}_{\delta}$.

\begin{theorem} \label{energydiscrete} Let $A \subset {\mathbb R}^d$ be a $1$-separated set of cardinality $N$. Let $\delta < \frac{1}{4 diam(A)}$, and let
\begin{equation}\label{ProbabilityMeasureMu}
d \mu(x)=N^{-1} \delta^{-d} \sum_{a \in A}
\chi_{B} \left(\delta^{-1} \left( x-\frac{a}{diam(A)} \right) \right) dx,
\end{equation}
where $\chi_{B}$ denotes the characteristic function of the ball of radius one centered at the origin.

Then
$$ I_{\alpha}(\mu)=I+II,$$ where
$$ I \approx N^{-1} \delta^{-\alpha},$$ and
$$ II \approx {(diam(A))}^{\alpha} \cdot N^{-2} \sum_{a \not=a'} {|a-a'|}^{-\alpha}.$$
Notice that the sum in $II$ is actually a double sum, in $a$ and $a'$.
\end{theorem}


\begin{proof}

By $B(x,r)$ we denote, as usual, the Euclidean ball of center $x$ and radius $r$. Then we split the energy integral in the diagonal and off-diagonal terms as follows

\begin{eqnarray}
I_{\alpha}(\mu) & = & \frac{1}{N^2 \; \delta^{2d}}  \sum_{a,a' \in A}  \int \int \frac{1}{{|x-y|}^{\alpha}} \; \chi_{B(\frac{a}{diam(A)},\delta )} (x) \; \chi_{B(\frac{a'}{diam(A)},\delta )} (y) \; dx \; dy = \nonumber \\
 & = & \sum_{a \in A} + \sum_{a \neq a' } = I + II  \nonumber
\end{eqnarray}

And direct calculations and estimates show that

$$
I \approx \frac{1}{N^2 \; \delta^{2d}} \; \delta^d  \left( \int_0^\delta  \frac{r^{d-1}}{r^\alpha} \; dr  \right) \; N \approx \frac{1}{N \; \delta^\alpha}
$$

and that

$$
II \approx \frac{1}{N^2 \; \delta^{2d}} \; \sum_{a \neq a' } \frac{ (diam(A))^\alpha }{{|a-a'|}^{\alpha}} \; \delta^{d} \;\delta^{d} \approx {(diam(A))}^{\alpha} \cdot N^{-2} \sum_{a \not=a'} {|a-a'|}^{-\alpha}
$$

\end{proof}

This leads us to a definition of Hausdorff $\alpha$-adaptability.

\begin{definition} \label{DefinitionHausdorffAlphaAdaptability}
Let $A_N \subset {\mathbb R}^d$ be a family of $1$-separated sets in ${\mathbb R}^d$, so that the cardinality of $A_N = \#\{A_N \} =N $. Assume that $${({[diam(A_N)]}^{-1}A_N)}_{\frac{1}{4diam(A_N)}} \subset {[0,1]}^d.$$
We say that $A_N$ is Hausdorff $\alpha$-adaptable if (\ref{diameterlim}) holds, that is $\delta \gtrsim N^{-\frac{1}{\alpha}}$ (with constant independent of $N$), and
\begin{equation}\label{OffDiagonalTermInEnergyIntegral}
{\cal I}_{\beta}(A_N)=N^{-2} \sum_{a \not=a'} {|a-a'|}^{-\beta} \lesssim {(diam(A_N))}^{-\beta},
\end{equation}
(also with constant independent of $N$, but that could depend on $\beta$), for all $\beta < \alpha.$
\end{definition}

\vskip.125in

Notice that the inequality $\geq$ always holds in \eqref{OffDiagonalTermInEnergyIntegral}. What \eqref{OffDiagonalTermInEnergyIntegral} says is that the average of the summands is actually comparable to the smallest summand.

The requirement that \eqref{OffDiagonalTermInEnergyIntegral} holds for all $\beta < \alpha$ is consistent with the continuous case where, although there is only one Hausdorff dimension for a set, call it $\alpha_0$, for any $0< \alpha < \alpha_0$, there exists a measure $\mu$ so that the energy integral  $I_\alpha(\mu)<\infty $ (this is a consequence of Frostman's lemma, see Theorem \ref{EnergyCharacterizationOfHausdorffDimension}.)

Although it is not part of Definition \ref{DefinitionHausdorffAlphaAdaptability}, later in the paper we will occasionally also work with the condition
\begin{equation}\label{OffDiagonalTermInEnergyIntegralWithLessapprox}
{\cal I}_{\alpha}(A_N)=N^{-2} \sum_{a \not=a'} {|a-a'|}^{-\alpha} \lessapprox {(diam(A_N))}^{-\alpha}.
\end{equation}

Notice that condition $\delta \gtrsim N^{-\frac{1}{\alpha}}$ is indeed condition \eqref{diameterlim}. Indeed, if given a set $A$ of cardinality $N$ we first rescale it by $\frac{1}{diam(A)}$, and then impose the condition that $\delta \leq \frac{1}{4 diam(A)}$, as summarized in the expression for $\mu$ in equation \eqref{ProbabilityMeasureMu}, then $\delta \approx \frac{1}{diam(A)}$, and \eqref{diameterlim} is equivalent to the condition $\delta \gtrsim N^{-\frac{1}{\alpha}}$, which is equivalent to saying that the diagonal term $I$ in Theorem \ref{energydiscrete} is bounded.

As with Minkowski dimension, when we allow for $\alpha$-adaptability of large subsets, we get the definition of Hausdorff dimension.

\begin{definition} \label{DiscreteHausdorffDimension}
We define the discrete Hausdorff dimension of a family of $1$-separated sets $A_N \subset {\mathbb R}^d$ with $\#\{A_N \} =N $ to be
\begin{eqnarray}
\dim_{\mathcal{H}} (A_N) = \sup \{ \beta >0 & : & \text{ for every } \varepsilon>0, \text{ there exists a family of sets } B_N \subseteq A_N \nonumber \\
& &  \text{ and a constant } C_\varepsilon >0, \text{ so that } \#(B_N) \geq \frac{C_\varepsilon}{N^\varepsilon} \#(A_N), \nonumber \\
& &  \text{ and so that } B_N  \text{ is Hausdorff $\beta$-adaptable. } \} \nonumber
\end{eqnarray}
The constant $C_\varepsilon$ depends on $\varepsilon$ and on the sequence $\{B_N \}$, but not on $N$ (and hence, the constant in \eqref{OffDiagonalTermInEnergyIntegral} ends up depending on $\varepsilon$ and on $\beta$ but not on $N$ when we compute the discrete Hausdorff dimension, since we have to check \eqref{OffDiagonalTermInEnergyIntegral} for all the possible $B_N$.) If there are no such $\beta >0$, the Hausdorff dimension of $A_N$ is zero.
\end{definition}

\vskip.125in

Notice also that if the condition $\delta \gtrsim N^{-\frac{1}{\alpha}}$ is satisfied for a certain $\alpha_0 >0$, then it is satisfied for all $0< \alpha < \alpha_0$ (see Theorem \ref{AlphaAdaptableIsMoreRestrictiveAsAlphaIncreases} below.) As a consequence, among the possible values of $\alpha$ for which the diagonal term $I$ in Theorem \ref{energydiscrete} is bounded, when looking for the $\alpha$ for which $A_N$ is Hausdorff $\alpha$-adaptable (if it exists), we look for the $\alpha$ that makes the off-diagonal term $II$ in Theorem \ref{energydiscrete} bounded. Considering these observations for all possible families of ``large subsets" $B_N$, we get that also in the discrete setting, $\dim_{\mathcal{H}} (A_N) \leq  \dim_{\mathcal{M}} (A_N)$. (It is in order to get this property that, given that we wanted to allow for ``large subsets" $B_N$ in the definition of discrete Hausdorff dimension, we also allowed for them in the definition of discrete Minkowski dimension.)


\begin{theorem}\label{AlphaAdaptableIsMoreRestrictiveAsAlphaIncreases}
Let $A_N \subset \mathbb{R}^d$ be a family of $1$-separated sets in $\mathbb{R}^d$, so that $\#{A_N}=N$. If $A_N$ is adaptable to the discrete Minkowski dimension $\alpha_0$, then $A_N$ is adaptable to the discrete Minkowski dimension $\alpha$, for any $0 \leq \alpha < \alpha_0$. If $A_N$ is Hausdorff $\alpha_0$-adaptable, then $A_N$ is Hausdorff $\alpha$-adaptable, for any $0 \leq \alpha < \alpha_0$.
\end{theorem}

\begin{proof}
Condition \eqref{diameterlim} is equivalent, as we have seen, to $\delta \gtrsim N^{- \frac{1}{\alpha}}$, for $\delta$ the minimum separation between two points in $A_N$, after $A_N$ has been rescaled to have diameter $\approx 1$. Notice now that $\alpha \to N^{- \frac{1}{\alpha}}$ is an increasing function of $\alpha$.

Notice also that $\displaystyle{ II \approx  {(diam(A))}^{\alpha} \cdot N^{-2} \sum_{a \not=a'} {|a-a'|}^{-\alpha} = \frac{1}{N^2} \sum_{a \not=a'} \left( \frac{ \max |a-a'| }{|a-a'| } \right)^\alpha }$, (in Theorem \ref{energydiscrete}), and that for $b>1$, the function $x \to b^x$ is increasing and positive, hence so is the last term in the previous equation.

\end{proof}

Our next Theorem is also related to the statement in the continuous case that for a set $E \subset \mathbb{R}^d$, $\dim_{{\cal H}}(E) \leq \dim_{{\cal M}}(E)$. In the sense that, although we already know that in the discrete setting, the Minkowski dimension is larger than the Hausdorff dimension, it might look as if this is so only because of the ``artificial" constraint of imposing condition \eqref{diameterlim} as part of Definition \ref{DefinitionHausdorffAlphaAdaptability}. Theorem \ref{MinkowskiDimensionGreaterThanHausdorffDimension} below shows that it is not such an ``artificial" requirement.


\begin{theorem}\label{MinkowskiDimensionGreaterThanHausdorffDimension}
Let $A_N \subset \mathbb{R}^d$ be a family of $1$-separated sets in $\mathbb{R}^d$, so that $\#{A_N}=N$. If equation \eqref{OffDiagonalTermInEnergyIntegral} is satisfied for a given $\alpha > 0$, (i.e. the control of the off-diagonal term in the energy integral), then equation \eqref{diameterlim} is satisfied in the same sense for the same $\alpha > 0$, at least by a subset of $A_N$ of size $\frac{N}{2}$ (i.e. the control of the diagonal term in the energy integral, or equivalently, the Minkowski dimension estimate.) An analogous statement holds with condition \eqref{OffDiagonalTermInEnergyIntegralWithLessapprox} instead of condition \eqref{OffDiagonalTermInEnergyIntegral}.

More precisely,
\begin{enumerate}
\item[(a)] If
$$
{\cal I}_{\alpha}(A_N)=\frac{1}{N^2} \sum_{a \not=a'} {|a-a'|}^{-\alpha} \lesssim 1,
$$
then, after rescaling to the unit cube in $\mathbb{R}^d$, and perhaps removing a subset of size at most $\frac{N}{2}$, the minimum separation between points $\delta$ satisfies $\delta \gtrsim N^{- \frac{1}{\alpha}}$.
\item[(b)] If
$$
{\cal I}_{\alpha}(A_N)=\frac{1}{N^2} \sum_{a \not=a'} {|a-a'|}^{-\alpha} \lessapprox 1,
$$
then, after rescaling to the unit cube in $\mathbb{R}^d$, and perhaps removing a subset of size at most $\frac{N}{2}$, the minimum separation between points $\delta$ satisfies $\delta \gtrapprox N^{- \frac{1}{\alpha}}$.
\end{enumerate}

\end{theorem}

\begin{proof}
First rescale $A_N$ to have diameter $1$. Then, in order to prove case $(a)$, we (essentially) want to prove that if

\begin{equation}\label{ControlOffDiagonalEquationRescaledToDiameter1}
{\cal I}_{\alpha}(A_N)=\frac{1}{N^2} \sum_{a \not=a'} {|a-a'|}^{-\alpha} \lesssim 1,
\end{equation}

then the minimum separation between points $\delta$ satisfies $\delta \gtrsim N^{- \frac{1}{\alpha}}$.

Notice first that if \eqref{ControlOffDiagonalEquationRescaledToDiameter1} is satisfied by $A_N$, then it is also satisfied (with slightly different constants) by any subset $B \subset A_N$ with $\#(B) \geq \frac{N}{2}$ (but the constants are the same for all such $B$.) So, let us fix a small $\varepsilon >0$, and assume it is not true that $\delta \geq \varepsilon N^{- \frac{1}{\alpha}}$ for $A_N$. Then there exist $a, a' \in A_N$ such that $| a-a' | \leq \varepsilon N^{- \frac{1}{\alpha}}$. Remove $a'$ from $A_N$, let the resulting set be $B_1$, and let us say that $a'$ no longer relates to $a$. If $B_1$ satisfies $\delta \gtrsim (N-1)^{- \frac{1}{\alpha}}$, stop since we are done. Otherwise, by the same reasoning, remove another point from $B_1$ thus yielding the set $B_2$. Continue in this manner for $\frac{N}{2}$ steps. If we have stopped at or before $\frac{N}{2}$ steps, we are done. If that is not the case, then, if we denote $E = \{ (a,a'): a' \text{ no longer relates to }  a  \}$, so that $\#(E) = \frac{N}{2}$, then going back to the original set $A_N$,

$$
\frac{1}{N^2} \sum_{a \not=a'} {|a-a'|}^{-\alpha} \geq \frac{1}{N^2} \sum_{(a,a') \in E} {|a-a'|}^{-\alpha} \geq \frac{1}{N^2} \frac{N}{2} \frac{N}{\varepsilon^\alpha} = \frac{1}{2 \varepsilon^\alpha}.
$$

Now letting $\varepsilon \to 0$, gives the desired contradiction.

The proof for case $(b)$ is completely analogous.

\end{proof}


\vskip.125in

\section{Applications of $\alpha$-adaptability to the Erd\H os-Falconer distance problem}\label{ApplicationsToErdosFalconer}

\vskip.125in

As we mentioned in the Introduction, the Erd\H os distance conjecture in geometric combinatorics says that if $A \subset {\mathbb R}^d$, $d \ge 2$, then
$$ \# \Delta(A) \gtrapprox {(\# A)}^{\frac{2}{d}}, $$ where
$$ \Delta(A)=\{|a-a'|: a,a' \in A\},$$ with
$$ {|x|}^2=x_1^2+x_2^2+\dots+x_d^2.$$

Taking $A={[0,N^{\frac{1}{d}}]}^d \cap {\mathbb Z}^d$ shows that one cannot in general do better. In the continuous setting, the Falconer distance conjecture says that if the Hausdorff dimension of
$E \subset {[0,1]}^d$ is larger than $\frac{d}{2}$, then the Lebesgue measure of $\Delta(E)$ is positive. Once again taking $A$ to be a set built on an appropriately scaled version of the integer lattice shows that it is possible for $\Delta(E)$ to have Lebesgue measure $0$ if the Hausdorff dimension of $E$ is any number less than [or equal for $d=2$] $\frac{d}{2}$. The connection between the putative sharpness examples in the two problems eventually led the first named author and I. Laba \cite{iosevichlabakdistance} to establish the following result in the special case of Delone sets.

\begin{definition} \label{delone} We say that $A \subset {\mathbb R}^d$ is Delone if there exist $C,c>0$ such that $A$ is $c$-separated and every cube of side-length $C$ contains at least one point of $A$. \end{definition}

For the purposes of this paper, we may prune and scale $A$ such that for every $m \in {\mathbb Z}^d$, $m+{[0,1]}^d$ contains exactly one point of $A$.

\begin{theorem} \label{iosevichlaba} \cite{iosevichlabakdistance} Let $A$ be a Delone set and define $A_q=A \cap {[0,q]}^d$. Suppose that the Falconer distance conjecture holds to the extent that if the Hausdorff dimension of $E \subset {[0,1]}^d$ is greater than $s_0$ ($s_0 \ge \frac{d}{2}$), then the Lebesgue measure of
$\Delta(E)$ is positive. Then
$$ \# \Delta(A_q) \gtrapprox q^{\frac{d}{s_0}}.$$
\end{theorem}

In particular, if $s_0=\frac{d}{2}$, as conjectured, then we see that the Falconer conjecture implies the Erd\H os conjecture in the context of Delone sets.

Let us now prove that Delone sets are Hausdorff $d$-adaptable, i.e. Theorem \ref{DeloneSetsAreDAdaptableVersionIntroduction}.

\begin{theorem} \label{DeloneSetsAreDAdaptable}
Delone sets have discrete Hausdorff dimension $d$ in $\mathbb{R}^d$.
\end{theorem}

\begin{proof}
Let $A$ be a Delone set in $\mathbb {R}^d$, and rescale it so that it is $1$-separated. Consider $A_N = A \cap [0,L_N]^d$ so that $\#(A_N) =N$. Then $L_N \approx  N^{\frac{1}{d}}$, since every cube of sidelength $C$ contains at least one point of $A$. Consequently, $diam(A_N) \lesssim  N^{\frac{1}{d}}$, which is condition \eqref{diameterlim}.

Notice that condition \eqref{OffDiagonalTermInEnergyIntegral} is scale invariant. Then, since $A_N$ is $1$-separated, and since each point $a \in A$ contributes the same amount to ${\cal I}_{\alpha}(A_N)$, up to comparability constants, and that amount can be calculated, again up to comparability constants by an integral which is computed by changing to polar coordinates, we get that for $0<\alpha<d$,

$$
N^{-2} \sum_{a \not=a'} {|a-a'|}^{-\alpha} \approx  \frac{1}{ N^{2}} \; N \int_1^{L_N} r^{d-1-\alpha } \; dr \approx
\frac{1}{N} \left( N^{\frac{1}{d}}  \right)^{d-\alpha} = N^{ \frac{-\alpha}{d}} \approx [diam(A_N)]^{-\alpha}
$$

hence $A$ (or $A_N$) is Hausdorff $\alpha$-adaptable, for $0<\alpha \leq d$, and thus, $\dim_{\mathcal{H}}(A)=d.$

As a curiosity, notice that for $\alpha=d$, we already know that condition \eqref{diameterlim} is satisfied, but also condition \eqref{OffDiagonalTermInEnergyIntegralWithLessapprox} is satisfied:


$$
N^{-2} \sum_{a \not=a'} {|a-a'|}^{-d} \approx \frac{1}{ N^{2}} \; N \int_1^{L_N} r^{-1}  \; dr \approx
\frac{1}{N} \left[ \log (N) \right] \lessapprox  \frac{1}{N}  \approx {(diam(A_N))}^{-d}
$$

\end{proof}

Notice that, for a $1$-separated set $A_N \subset \mathbb{R}^d$ of cardinality $N$, the minimum diameter of $A$ among such sets, is precisely comparable to $N^{\frac{1}{d}}$ (attained when all points are packed roughly in a lattice, i.e. precisely in the case of a Delone set.) This simple remark proves that the discrete Hausdorff dimension (and Minkowski dimension) of such a set is always $\leq d$, as in the continuous case.

We now prove Theorem \ref{HausdorffAlphaAdaptableSatisfyErdosIfAlphaAccordingToFalconerVersion2WhenSetsAreNestedVersionIntroduction}.

\begin{theorem} \label{HausdorffAlphaAdaptableSatisfyErdosIfAlphaAccordingToFalconerVersion2WhenSetsAreNested}
Suppose that the Falconer distance conjecture holds to the extent that if the Hausdorff dimension of $E \subset {[0,1]}^d$ is greater than $s_0$ ($s_0 \ge \frac{d}{2}$), then the Lebesgue measure of $\Delta(E)$ is positive.

\begin{enumerate}
\item[(a)] Let $A_N$ be a family of sets with $\# (A_N) = N$ which is Hausdorff $\alpha_0$-adaptable, for some $\alpha_0 > s_0$. 
Assume also that for any $s_0 < \alpha < \alpha_0$, the family $ C_{N,\alpha}={({[diam(A_N)]}^{-1}A_N)}_{N^{-\frac{1}{\alpha}}}$ is a nested family of sets, i.e. $C_{N+1,\alpha} \subseteq C_{N,\alpha}$.
Then
$$ \# \Delta(A_N) \gtrapprox N^{\frac{1}{s_0}}.$$

\item[(b)] Let $A_N$ be a family of sets with $\# (A_N) = N$. Assume also that for any $\alpha$ with $s_0 < \alpha < \alpha_0$, and for every $\widetilde{\varepsilon} >0$ there exists a family of subsets $B_N \subseteq A_N$ and a constant $C_{\widetilde{\varepsilon}} >0$ (which depends on $\widetilde{\varepsilon}$, on $\alpha$, and on the sequence $\{ B_N \}$, but not on $N$), so that
$\# (B_N) \geq \frac{ C_{ \widetilde{\varepsilon} } }{ N^{ \widetilde{\varepsilon}} } \# (A_N)$
, and $B_N$ is Hausdorff $\alpha$-adaptable, and the family
$ C_{N,\alpha}={({[diam(B_N)]}^{-1}B_N)}_{(\#(B_N))^{-\frac{1}{\alpha}}}$ is a nested family of sets, i.e. $C_{N+1,\alpha} \subseteq C_{N,\alpha}$.
Then
$$ \# \Delta(A_N) \gtrapprox N^{\frac{1}{s_0}}.$$

\end{enumerate}

\end{theorem}


\begin{proof}

Let us first prove part $(a)$. Let us assume, for a contradiction, that $\# \left( \Delta (A_N) \right) $ is not $ \gtrapprox N^{ \frac{1}{s_0}}$, i.e. that there exists an $\varepsilon >0$ and a subsequence $A_{N_j}$ with
\begin{equation}\label{SmallNumberOfDistances}
\# \left( \Delta (A_{N_j}) \right) < N_j^{\frac{1}{s_0} - \varepsilon  }.
\end{equation}


Take now an $\alpha > s_0$ but so close to $s_0$ that $\frac{1}{s_0} -\varepsilon < \frac{1}{\alpha}$ (which we can do by Theorem \ref{AlphaAdaptableIsMoreRestrictiveAsAlphaIncreases}.) Recall now from \eqref{ProbabilityMeasureMu} that, associated to each $A_N$, we have the probability measure

\begin{equation}\label{ProbabilityMeasureMuForAN}
d \mu_{A_N}(x)= \frac{c}{N} \delta^{-d} \sum_{a \in A}
\chi_{B} \left(\delta^{-1} \left( x-\frac{a}{diam(A)} \right) \right) dx,
\end{equation}
where $\chi_{B}$ denotes the characteristic function of the ball of radius one centered at the origin, and $c$ is an absolute constant that does not depend on $N$ (it actually only depends on the volume of the unit ball in $\mathbb{R}^d$.) We pick $\delta \approx N^{- \; \frac{1}{\alpha}}$.

If we call the support of $\mu_{A_N}$, $supp (\mu_{A_N}) = K_N \subset [-1,2]^d$, by the Blaschke selection theorem (see e.g. \cite{falconerbook} p.37), there is a further subsequence of the $K_{N_j}$, which we will keep calling $K_{N_j}$ for simplicity, so that $K_{N_j} \to \widetilde{K_0}$, with convergence in the Hausdorff metric.
There is a further subsequence of the family of sets $A_{N_j}$, which again we keep calling $A_{N_j}$, so that the measures $\mu_{A_{N_j}}$ converge weakly (using the measure-theoretic terminology, in functional analysis the term would be weak-$\ast$ convergent). So we have that $\mu_{A_{N_j}} \rightharpoonup \mu_0$.

Then we claim that

\begin{equation}\label{SupportWeakLimitMeasureAndHausdorffLimitSupports}
K_0 := supp (\mu_0) \subseteq \widetilde{K_0},
\end{equation}

although equality need not hold. In order to prove \eqref{SupportWeakLimitMeasureAndHausdorffLimitSupports}, let $x_0 \in supp( \mu_0)$. Then, for every $\eta >0$, $\mu_0 (B(x_0,\eta ))>0$, where $B(x_0,\eta )$ denotes the open ball of center $x_0$ and radius $\eta$. Then (see e.g. \cite{mattila} p.19),

$$
\liminf_{N_j \to \infty}  \mu_{A_{N_j}} ( B(x_0,\eta ) ) \geq \mu_0 ( B(x_0,\eta ) ) >0 ,
$$

so for any $N_j$ sufficiently large, there is a point $a_{N_j,x_0} \in A_{N_j} \cap \overline{B \left(x_0, \eta + N_j^{- \;\frac{1}{\alpha}} \right)}$. Taking $\eta  \to 0$ and $N_j \to \infty$, we have that $supp (\mu_{A_{N_j}}) \supset A_{N_j} \ni a_{N_j,x_0} \to x_0$, and hence $x_0 \in \widetilde{K_0}$.

On the other hand, since the family $A_N$ is Hausdorff $\alpha$-adaptable,
by Theorem \ref{energydiscrete}, the energy integrals $I_\alpha (\mu_{A_{N_j}}) \leq C < \infty$ (with $C$ independent of $N_j$.) A well-known lemma in potential theory then yields that

\begin{equation}\label{LimitMeasureHasFiniteEnergyIntegral}
I_\alpha (\mu_0) \leq C < \infty.
\end{equation}

For the convenience of the reader, we now sketch the main ideas in the proof of the aforementioned lemma. If $\mu_m \rightharpoonup \mu_0$, then $\mu_m \times \mu_m \rightharpoonup \mu_0 \times \mu_0$ (a consequence of the Stone-Weierstrass theorem). Use $\mu_m \times \mu_m \rightharpoonup \mu_0 \times \mu_0$ for each one of the continuous kernels $k_{\alpha,n} (x,y) = \min \left\{ \frac{1}{|x-y|^\alpha},n \right\}$, and apply the monotone convergence theorem.

As a consequence of \eqref{LimitMeasureHasFiniteEnergyIntegral} and Theorem \ref{EnergyCharacterizationOfHausdorffDimension}, recalling $K_0 := supp (\mu_0)$, then we have that $\dim_{\mathcal{H}} (K_0) \geq \alpha > s_0 \geq \frac{d}{2}$. Hence, Falconer's conjecture implies that

\begin{equation} \label{DistanceSetOfSupportOfLimitMeasureHasPostiveLength}
\mathcal{L}^1 (\Delta (K_0)) >0.
\end{equation}

Recalling $K_{N_j} := supp (\mu_{A_{N_j}}) $, it follows from the fact that $K_{N_j} \to \widetilde{K_0}$ in the Hausdorff metric, that $\Delta ( K_{N_j}  ) \to \Delta ( \widetilde{K_0} )$ in the Hausdorff metric. To see this, note that if $F_N \to F$ in the Hausdorff metric, then for every $\delta >0$, for a sufficiently large $N$, we have that $(F_N)_{\delta} \supseteq F$ and that $(F)_{\delta} \supseteq F_N$, so the same relations hold when taking $\Delta$. Now note that $\Delta (A_\delta) = \left(  \Delta (A)  \right)_{2\delta} $ .

Recall now that $\alpha > s_0$ was taken so close to $s_0$ that $\frac{1}{s_0} -\varepsilon < \frac{1}{\alpha}$. Due to the nesting of $(A_{N_j})_{\delta_j}$, where $\delta_j = N_j^{- \; \frac{1}{\alpha}}$, we have that  $ \left( \Delta (A_{N_j}) \right)_{2\delta_j} \supseteq \Delta (\widetilde{K_0})$, and then

\begin{equation}\label{InequalitiesProvingDistanceSetK0HasLengthZero}
\mathcal{L}^1 \left\{ \left( \Delta (A_{N_j}) \right)_{2\delta_j}   \right\} \geq \mathcal{L}^1 \left\{ \Delta (\widetilde{K_0})  \right\} \geq  \mathcal{L}^1 \left\{ \Delta (K_0)  \right\}
\end{equation}

but $\mathcal{L}^1 \left\{ \left( \Delta (A_{N_j}) \right)_{2\delta_j}   \right\} \lesssim N_j^{\frac{1}{s_0} -\varepsilon } \cdot N_j^{- \; \frac{1}{\alpha}} \to 0 $, which proves that $\mathcal{L}^1 \left\{ \Delta (K_0)  \right\} =0$, a contradiction with \eqref{DistanceSetOfSupportOfLimitMeasureHasPostiveLength}.

With respect to part $(b)$, let us remark that because of the nesting property of the family $B_N$, the statement is assuming something actually stronger than saying that $\dim_{\mathcal{H}} (A_N) \geq \alpha_0$. The proof of part $(b)$ is the same as that of part $(a)$, only substituting $A_N$ for $B_N$, and $N$ for $C_{ \widetilde{\varepsilon} } N^{1-\widetilde{\varepsilon} }$ (analogously for $N_j$.) Then the proof of part $(a)$ yields
$$
\# (A_N) \geq \# (B_N) \geq C_{ \widetilde{\varepsilon} } N^{( 1-\widetilde{\varepsilon} ) \frac{1}{s_0} }
$$
and since this holds for every $\widetilde{\varepsilon} >0$, the result follows.

\end{proof}

\begin{remark}
As a curiosity, in order to see that equality need not hold in \eqref{SupportWeakLimitMeasureAndHausdorffLimitSupports}, take $M$ points uniformly distributed in $[0,1] \times [\frac{1}{2},1]$, and take $M^2$ points uniformly distributed in $[0,1] \times \{0\}$. Let $N = M + M^2$ and let $A_N$ be the union of those points. Then it is easy to see that the points on $[0,1] \times \{0\}$ outweigh substantially the points in $[0,1] \times [\frac{1}{2},1]$, to the point that for any weakly convergent subsequence $\mu_{A_{N_j}} \rightharpoonup \mu_0$, we have that $supp (\mu_0) = [0,1] \times \{0\} \subsetneq [0,1] \times [\frac{1}{2},1] = \widetilde{K_0}$. This curiosity highlights the fact that in the machinery being developed in this paper, it is important not only what set the sequence $A_N$ approaches, but also how it approaches this set, in the sense of with what weights it approaches it.

\end{remark}

We also get another ``translation theorem" from Falconer to Erd\H os, without the assumption that the sets are nested, but with an extra assumption in the form of a slightly stronger version of the Falconer conjecture, namely that not only the distance set $\Delta(E)$ has positive length, but that there is a quantitative control of the length $\mathcal{L}^1 (\Delta(E)) \geq C=C(\alpha, C_0)>0$ (see below for the meaning of these parameters.) However, as we noted in the Introduction, all known recent proofs of results pertaining to the Falconer conjecture actually yield such a quantitative control of the length. We prove it in a slightly more general form than Theorem \ref{HausdorffAlphaAdaptableSatisfyErdosIfAlphaAccordingToFalconerVersion2WhenFalconerIsStrengthenedVersionIntroduction}.

\begin{theorem} \label{HausdorffAlphaAdaptableSatisfyErdosIfAlphaAccordingToFalconerVersion2WhenFalconerIsStrengthened}

\begin{enumerate}
\item[(a)]
Suppose that the Falconer distance conjecture holds to the extent that if a Borel probability measure $\mu$ supported on $E \subset {[0,1]}^d$ satisfies that $I_\alpha (\mu) \leq C_0 < \infty$, for some $\alpha > s_0 \ge \frac{d}{2}$ (recall Theorem \ref{EnergyCharacterizationOfHausdorffDimension}), then $\mathcal{L}^1 (\Delta(E)) \geq C=C(\alpha, C_0)>0$.

Let $A_N \subset \mathbb{R}^d$ be a family of sets with $\# (A_N) = N$ with $\dim_{\mathcal{H}}(A_N) = \alpha_0 > s_0$.
%
Then
$$ \# \Delta(A_N) \gtrapprox N^{\frac{1}{s_0}}.$$

(Slightly) more generally, let $A_N \subset \mathbb{R}^d$ be a family of sets with $\# (A_N) = N$ such that, for every $\widetilde{\varepsilon} >0$, there exists a family of subsets $B_N \subseteq A_N$ and a constant $C_{ \widetilde{\varepsilon} }$ (which may depend on $\widetilde{\varepsilon}$, and the sequence $\{B_N\}$, but not on $N$), with $\#(B_N) \geq \frac{ C_{ \widetilde{\varepsilon} } }{ N^{ \widetilde{\varepsilon} } } \#(A_N)$, so that $B_N$ satisfies equation \eqref{OffDiagonalTermInEnergyIntegral} for some $\alpha_0 > s_0$ (with constant that may depend on $\alpha_0$, $\widetilde{\varepsilon}$, and the sequence $\{B_N\}$, but not on $N$.) Then
$$ \# \Delta(A_N) \gtrapprox N^{\frac{1}{s_0}}.$$

\item[(b)]
Assume the Falconer distance conjecture holds to the extent that for any Borel probability measure $\mu$ supported on $E \subset {[0,1]}^d$ that satisfies that $I_\alpha (\mu) \lessapprox 1$, for some $\alpha > s_0 \ge \frac{d}{2}$, then $\mathcal{L}^1 (\Delta(E)) \geq C=C(\alpha, C_0)>0$.

Let $A_N \subset \mathbb{R}^d$ be a family of sets with $\# (A_N) = N$ such that, for every $\widetilde{\varepsilon} >0$, there exists a family of subsets $B_N \subseteq A_N$ and a constant $C_{ \widetilde{\varepsilon} }$ (which may depend on $\widetilde{\varepsilon}$, and the sequence $\{B_N\}$, but not on $N$), with $\#(B_N) \geq \frac{ C_{ \widetilde{\varepsilon} } }{ N^{ \widetilde{\varepsilon} } } \#(A_N)$, so that $B_N$ satisfies equation \eqref{OffDiagonalTermInEnergyIntegralWithLessapprox} for some $\alpha_0 > s_0$ (with constant that may depend on $\alpha_0$, $\widetilde{\varepsilon}$, and the sequence $\{B_N\}$, but not on $N$.) Then
$$ \# \Delta(A_N) \gtrapprox N^{\frac{1}{s_0}}.$$
\end{enumerate}

\end{theorem}

\begin{proof}
Fix $\widetilde{\varepsilon} >0$. Regarding part $(a)$, with the same notation as in \eqref{ProbabilityMeasureMuForAN}, by Theorems \ref{AlphaAdaptableIsMoreRestrictiveAsAlphaIncreases} and \ref{MinkowskiDimensionGreaterThanHausdorffDimension}, if necessary after removing a subset of size at most $\frac{\#(B_N)}{2}$ from $B_N$ (but we will keep calling the resulting set $B_N$), we get $I_{\alpha} (\mu_{B_N}) \leq C' \cdot C_0 < \infty$ for any $\alpha \leq \alpha_0$ (where $C'$ is an absolute constant.) Hence, for $\delta_\alpha = (\#(B_N))^{- \; \frac{1}{\alpha}}$, we have that $\mathcal{L}^1 (\Delta( (B_N)_{\delta_\alpha}  )) \geq C=C(\alpha, C_0)>0$.

Then the number of different Euclidean distances determined by $A_N$ satisfies
$$\# \Delta(A_N) \geq \# \Delta(B_N) \gtrsim \frac{C}{\delta_\alpha} = C \left[ \#(B_N) \right]^{-\frac{1}{\alpha}} \gtrsim  N^{ \left( 1-  \widetilde{\varepsilon} \right) \left( \frac{1}{s_0} - \varepsilon  \right) },$$
for any $\varepsilon >0$ (by taking $\alpha$ as close as we want to $s_0$.) Now send both $\varepsilon$ and $\widetilde{\varepsilon}$ to zero.

The proof for part $(b)$ is analogous.

\end{proof}

\section{Examples}\label{Examples}

Our next Theorem (mentioned in the Introduction as Theorem \ref{FeketePointsApproximationVersionIntroduction}) shows that there are plenty of cases to which our machinery applies (and also plenty of them to which it does not apply, at least directly, in the sense that a priori it is possible to find a ``sufficiently large" subset inside the following examples to which our machinery could be applied to calculate distances, as in the example from Theorem \ref{ExampleMinkowskiDimensionDifferentToHausdorffDimension} below.)

\begin{theorem}\label{FeketePointsApproximation}
Let $E \subset [0,1]^d$ be a compact set with diameter $diam (E) \approx 1$, so that there exists a Borel probability measure $\mu$ supported on $E$ with $I_\alpha (\mu) < \infty$ (see Theorem \ref{EnergyCharacterizationOfHausdorffDimension}), for $0< \alpha <d$. Then there exists a family of Hausdorff $\alpha$-adaptable sets $A_{N_j} \subset [0,1]^d $, with $\# (A_{N_j} ) = N_j \to \infty$, so that, with the notation of \eqref{ProbabilityMeasureMuForAN}, $\mu_{A_{N_j}} \rightharpoonup \mu_0$ (weak-$\ast$ convergence) with $\mu_0$ a Borel probability measure supported on $K_0$ satisfying $I_\alpha (\mu_0) < \infty$, and $A_{N_j} \to \widetilde{K_0}$ in the Hausdorff metric, with $K_0 \subseteq \widetilde{K_0}  \subseteq E$.
\end{theorem}

\begin{proof}
A possible approach to this Theorem is to discretize the construction of the Frostman measure. However, this Theorem is essentially already known in the literature as the Fekete-Szeg\H o theorem (see \cite{ransford}) or transfinite diameter (see also \cite{landkof}.)

For the convenience of the reader, we recall the construction of the transfinite diameter and the proof that it equals the Riesz capacity, following \cite{landkof}, since we will need some elements of it.

Let $C_{\alpha} (E) = \sup \{ I_\alpha (\mu)^{-1} : \mu \text{   is a Radon probability measure with } supp(\mu) \subseteq E \}$, denote the Riesz capacity of order $\alpha$ of $E$. From the hypotheses, $C_{\alpha} (E) >0$.

Consider the function
\begin{equation}\label{FunctionForNthApproximantOfTransfiniteDiameter}
F_{\alpha}(x_1, \dots , x_N) =   \frac{1}{  \binom{N}{2} } \sum_{i < j} \frac{1}{|x_i - x_j|^\alpha}
\end{equation}
defined on $E \times \dots \times E.$

Since $E$ is compact, $F_{\alpha}(x_1, \dots , x_N)$ achieves its minimum value on $E$ at certain points $x_i=\xi_i^{(N)}$. Let us define
\begin{equation}\label{NthApproximantOfTransfiniteDiameter}
D_N^{(\alpha)} =  \binom{N}{2} \left(  \sum_{i < j} \frac{1}{\left|\xi_i^{(N)} - \xi_j^{(N)}\right|^\alpha}  \right)^{-1}
\end{equation}

In order to compare the sum in $D_N^{(\alpha)}$ with $N$ elements and the $N$ possible sums for the subsets of $N-1$ elements, notice that

$$
\sum_{i < j} \frac{1}{\left|\xi_i^{(N)} - \xi_j^{(N)}\right|^\alpha} = \frac{1}{N-2} \sum_{k=1}^{N} \sum_{i < j}^{(k)} \frac{1}{\left|\xi_i^{(N)} - \xi_j^{(N)}\right|^\alpha}
$$

where $\displaystyle{\sum^{(k)}}$ denotes the sum in which the terms for $i=k$ and $j=k$ have been omitted. But

$$
\sum_{i < j}^{(k)} \frac{1}{\left|\xi_i^{(N)} - \xi_j^{(N)}\right|^\alpha} \geq \binom{N-1}{2} \frac{1}{D_{N-1}^{(\alpha)}},
$$

and consequently

$$
\frac{\binom{N}{2}}{D_N^{(\alpha)}} = \sum_{i < j} \frac{1}{\left|\xi_i^{(N)} - \xi_j^{(N)}\right|^\alpha}
\geq \frac{N}{N-2} \binom{N-1}{2} \frac{1}{D_{N-1}^{(\alpha)}} = \frac{\binom{N}{2}}{D_{N-1}^{(\alpha)}}.
$$

Therefore we get that
\begin{equation}\label{NthApproximantOfTransfiniteDiameterDecreases}
D_{N-1}^{(\alpha)} \geq  D_N^{(\alpha)},
\end{equation}

and hence $\displaystyle{ D^{(\alpha)} (E) := \lim_{N \to \infty} D_N^{(\alpha)} }$ exists (it is called the transfinite diameter of order $\alpha$ of $E$.)

Integrating the inequality
$$
\frac{\binom{N}{2}}{D_N^{(\alpha)}} \leq \sum_{i < j} \frac{1}{|x_i - x_j|^\alpha}
$$
against $d\nu (x_1) \dots d\nu (x_N)$, where $\nu$ is the equilibrium distribution on $E$ (in particular, by definition, a probability measure), gives

\begin{equation}\label{FirstInequalityTransfiniteDiameterCapacity}
D^{(\alpha)} (E) \geq C_\alpha (E).
\end{equation}

Consider the measure $\displaystyle{ \nu_N  = \frac{1}{N} \sum_{i=1}^{N} \delta_{\xi_i^{(N)}} }$, where $\delta_a$ is the Dirac delta measure at the point $a$.

This measure has infinite $\alpha$-energy $I_\alpha$, but if we use the truncated kernel
$$k_{\alpha,n} (x,y) = \min \left\{ \frac{1}{|x-y|^\alpha},n \right\}$$

then

\begin{equation}
\int_{E \times E} k_{\alpha,n} (x,y) d\nu_N (x) d\nu_N (y) \leq \frac{1}{N^2} \sum_{i \neq j} \frac{1}{\left|\xi_i^{(N)} - \xi_j^{(N)}\right|^\alpha} + \frac{n}{N} =
 \frac{2}{N^2}  \frac{\binom{N}{2}}{D_N^{(\alpha)}} + \frac{n}{N}
\end{equation}

Since $k_{\alpha,n} (x,y)$ is a continuous function, fixing $n$, by weak-$\ast$ compactness of measures, we may assume, passing to a subsequence, that $\nu_N \rightharpoonup \nu_0$. Then we obtain

\begin{equation}
\int_{E \times E} k_{\alpha,n} (x,y) d\nu_0 (x) d\nu_0 (y) \leq \frac{1}{D^{(\alpha)} (E)}
\end{equation}

Now applying the monotone convergence theorem gives $I_\alpha (\nu_0) \leq \frac{1}{D^{(\alpha)} (E)}$. Hence, using \eqref{FirstInequalityTransfiniteDiameterCapacity}, we get
$$
I_\alpha (\nu_0) \leq \frac{1}{D^{(\alpha)} (E)} \leq  \frac{1}{C_\alpha (E)} = I_\alpha (\nu),
$$
so that, by the uniqueness of the equilibrium distribution, $\nu_0 = \nu$ and

\begin{equation}\label{RelationTransfiniteDiameterCapacity}
D^{(\alpha)} (E) = C_\alpha (E).
\end{equation}

Consider now the family of sets $B_N = \left\{ \xi_i^{(N)} \right\}_{i=1}^{N}$ and the associated measures $\mu_{B_N}$, as in \eqref{ProbabilityMeasureMuForAN}. By the minimizing property of the $B_N$, we have that $diam (B_N) \approx 1$. If this were not the case, then $diam (B_N) << 1$, and by moving one of the points in $B_N$ as far as possible from the others (so that the diameter gets comparable to $1$), we would decrease the value in \eqref{FunctionForNthApproximantOfTransfiniteDiameter}. Notice that we are {\it not} stating that all points $ \xi_i^{(N)} \in \partial E $, where $\partial E$ is the boundary of $E$. This last statement is, in general, {\it false}. More precisely, if $\alpha > d-2$ in $\mathbb{R}^d$, the equilibrium distribution is in general {\it not concentrated on } $\partial E$ (see e.g. \cite{landkof} p.163.)

Since $diam (B_N) \approx 1$, by Theorem \ref{energydiscrete} and \eqref{NthApproximantOfTransfiniteDiameter}, the off-diagonal term $II$ in $I_\alpha (\mu_{B_N})$ is $\approx \frac{1}{D_N^{(\alpha)}}$, with absolute constants. By Theorem \ref{MinkowskiDimensionGreaterThanHausdorffDimension}, and again Theorem \ref{energydiscrete}, there exists a family of sets $A_N$ with $A_N \subseteq B_N$, and $\frac{N}{2} \leq \#(A_N) \leq N$, with $I_\alpha (\mu_{A_N}) \lesssim I_\alpha (\mu_{B_N}) $, again with absolute constants, since the sum in the term $II$ for $I_\alpha (\mu_{A_N}) $ has less terms than the corresponding sum for $\mu_{B_N}$.

By \eqref{NthApproximantOfTransfiniteDiameterDecreases} and \eqref{RelationTransfiniteDiameterCapacity}, $I_\alpha (\mu_{A_N}) \lesssim I_\alpha (\nu) = \frac{1}{C_\alpha (E)}$, again with absolute constants, so that the family $A_N$ is Hausdorff $\alpha$-adaptable. Note that the assumption $\alpha >0$ immediately implies that $\#(E) = \infty$.
By taking successive subsequences, we can assume that for a sequence of $N_j \to \infty$, $A_{N_j} \to \widetilde{K_0}$ in the Hausdorff metric, and $\mu_{A_{N_j}} \rightharpoonup \mu_0$ in weak-$\ast$ convergence. Then, as in \eqref{LimitMeasureHasFiniteEnergyIntegral}, $I_\alpha (\mu_0) < \infty $. If we call $K_0 = supp (\mu_0)$, then, as in \eqref{SupportWeakLimitMeasureAndHausdorffLimitSupports}, $K_0 \subseteq \widetilde{K_0}$. Also, since $A_{N_j} \subseteq B_{N_j} \subseteq E$, we have that $\widetilde{K_0} \subseteq E$.

\end{proof}

Our next Theorem gives an example of a family of sets $A_N \subset \mathbb{R}^d$ which is not Hausdorff $\alpha$-adaptable for any $\alpha>0$, and hence the machinery developed so far would seem not to apply at first sight in terms of producing Erd\H os type results assuming Falconer type results (if we had not introduced the considerations on large subsets of such families.) However, a closer look at the family of sets shows that the aforementioned machinery can indeed be applied, since indeed $\dim_{\mathcal{H}}(A_N)=d$.

\begin{theorem}\label{ExampleMinkowskiDimensionDifferentToHausdorffDimension}
There exists a family of $1$-separated sets $A_N \subset \mathbb{R}^2$, with $\#(A_N)=N$, which is Minkowski $1$-adaptable, but is not Hausdorff $\alpha$-adaptable, for any $\alpha >0$. However $\dim_{\mathcal{H}}(A_N)=2$ and hence, if the Falconer distance conjecture is true, then the family $A_N$ satisfies the Erd\H os distance conjecture $ \# \Delta(A_N) \gtrapprox N$, i.e. for any $\varepsilon >0$, there exists a constant $C_\varepsilon >0$, such that $$ \# \Delta(A_N) \geq C_\varepsilon N^{1-\varepsilon} .$$
\end{theorem}

\begin{proof}
For large $M$, let $B_M = \{\frac{1}{n}: n=1, \dots , M \}$, and let $A_N = B_M \times B_M$, with $N=M^2$. Rescale by $M^2$, so that the $x$ and $y$ coordinates of the points in the rescaled $A_N$ (let us call it $\widetilde{A_N}$) are precisely $M^2, \frac{M^2}{2}, \frac{M^2}{3}, \dots , \frac{M^2}{M-1}, M$. Then the minimum distance $\delta$ between two points in $\widetilde{A_N}$ is $\delta = \frac{M^2}{M-1}- M \approx 1$. Since $\diam(\widetilde{A_N}) = \sqrt{2}(M^2-M) \approx M^2 =N$, then $A_N$ is Minkowski $1$-adaptable.

Now, since equation \eqref{OffDiagonalTermInEnergyIntegral} is scale invariant, consider the interactions between  points of the form $a=(\frac{1}{p},\frac{1}{l}) \in A_N$ with points of the form $a'=(\frac{1}{n},\frac{1}{k}) \in A_N$, under the restrictions that $\frac{M}{10} \leq l,p \leq \frac{2M}{10}$, $n \geq \frac{M}{2}$, and $\frac{2M}{10} \leq k \leq \frac{3M}{10}$.

Consider the angle $\beta$ determined by $a'$, $a$, and the point $(0,\frac{1}{l})$. Then $0 \leq \beta \leq \beta_0$, where $\beta_0$ is the angle determined by $\left( \frac{1}{\frac{M}{2}}, \frac{1}{\frac{3M}{10}} \right)$, $\left( \frac{1}{\frac{2M}{10}}, \frac{1}{\frac{M}{10}} \right)$, and $\left( 0, \frac{1}{\frac{M}{10}} \right)$. Hence, $\tan (\beta_0) = \frac{20}{9}$, and for $0 \leq \beta \leq \beta_0$, $\cos (\beta) \geq \cos (\beta_0) \approx 0.41$, i.e. an absolute constant. Hence, if $P_{a,a'} = (\frac{1}{n},\frac{1}{l})$, we have that $| a-a' | \approx | a- P_{a,a'}|$ with universal constants that only depend on $\cos (\beta_0) \approx 0.41$.

As a consequence, if we fix $a$, and sum over all the described $a'$, since there are $\approx M$ possible values for $k$, and since $n>p>0$

$$
\sum_{a' \; : \; a \not=a'} \frac{1}{{|a-a'|}^\alpha} \approx M \sum_{ n \geq \frac{M}{2} } \frac{1}{{| \frac{1}{p} -\frac{1}{n}|}^\alpha} = M p^\alpha \sum_{ n \geq \frac{M}{2} } \frac{ n^\alpha }{( n-p)^\alpha } \geq M p^\alpha \frac{M}{2} \approx M^2 p^\alpha.
$$

If we now sum over $l$, but keeping $p$ fixed, since there are $\approx M$ such $l$, we get

$$
\sum_{l} \sum_{a' \; : \; a \not=a'} \frac{1}{{|a-a'|}^\alpha} \gtrsim M^3 p^\alpha.
$$

And now, summing over $p$,

$$
\sum_{a,a' \; : \; a \not=a'} \frac{1}{{|a-a'|}^\alpha} \gtrsim M^3 \sum_{p=\frac{M}{10}}^{\frac{2M}{10}}  p^\alpha \gtrsim M^{4+\alpha},
$$

since $\displaystyle{ \sum_{p=\frac{M}{10}}^{\frac{2M}{10}}  p^\alpha \approx \int_{\frac{M}{10}}^{\frac{2M}{10}} x^\alpha \; dx \approx M^{1+\alpha} }$.

Since $N=M^2$, and $\diam(A_N) \approx 1$, then for the whole set $A_N$ we have that

$$
II  \approx {(\diam(A_N))}^{\alpha} \cdot N^{-2} \sum_{a \not=a'} {|a-a'|}^{-\alpha} \gtrsim M^\alpha
$$

which is not bounded for any $\alpha >0$.

Although we do not need it, let us mention that a reasoning very similar to the one just done gives the upper bound $\displaystyle{ \sum_{\substack{ a,a' \in A_N \\ a \not=a' } } {|a-a'|}^{-\alpha} \lesssim M^{4+\alpha} }$, so that, indeed, $\displaystyle{ \sum_{\substack{ a,a' \in A_N \\ a \not=a' } } {|a-a'|}^{-\alpha} \approx M^{4+\alpha} }$. More precisely, consider $a=(\frac{1}{p},\frac{1}{l}) \in A_N$, and consider the lines that form an angle of $\frac{\pi}{4}$ with the coordinate axes through $a$, i.e., the lines $L_{a,1} \equiv x-y=\frac{1}{p} - \frac{1}{l}$, and $L_{a,2} \equiv x+y=\frac{1}{p} + \frac{1}{l}$. These lines divide the whole plane (and in particular the set $A_N$) into 4 sectors, denoted N,S,E,W (for North, South, East and West) in the obvious way. Let us consider a point $a'=(\frac{1}{n},\frac{1}{k}) \in A_N$ which is, say, in the W sector for $a$ (denoted $W(a)$). Define $P_a(a') = (\frac{1}{n},\frac{1}{l})$, i.e. the projection of $a'$ onto the line parallel to the coordinate axes in $W(a)$. Again by trigonometry, with universal constants, $|a-a'| \approx |a-P_a(a')|$. For a fixed $n$, there are at most $\approx M$ such points $a' \in W(a)$. The same reasoning applied to the other sectors for $a$ shows that for a fixed $a=(\frac{1}{p},\frac{1}{l}) \in A_N$, the interactions of $a$ with all other points $a'$ is bounded by $M$ times the interactions between $a$ and all other points $a'$ in the same row or column as $a=(\frac{1}{p},\frac{1}{l})$, i.e.

$$
\sum_{a' \; : \; a \not=a'} \frac{1}{{|a-a'|}^\alpha} \lesssim M  \left\{  \sum_{ \substack{a' \; : \; a \not=a' \\ a'=(\frac{1}{p},\frac{1}{k} ) } } \frac{1}{{|a-a'|}^\alpha} +  \sum_{ \substack{a' \; : \; a \not=a' \\ a'=(\frac{1}{n},\frac{1}{l} ) } } \frac{1}{{|a-a'|}^\alpha}   \right\}
$$

Let us focus on the interactions between $a=(\frac{1}{p},\frac{1}{l}) \in A_N$ and other points in its same row (the reasoning for the same column is symmetric.)

\begin{eqnarray}
\sum_{ \substack{ n \neq p \\ 1 \leq n \leq p } } \frac{1}{ |\frac{1}{p} - \frac{1}{n} |^\alpha } = p^\alpha \sum_{ \substack{ n \neq p \\ 1 \leq n \leq p } } \frac{ n^\alpha }{ |n-p|^\alpha } & = &  p^\alpha \left\{ \sum_{n=1}^{\frac{p}{2}-1} + \sum_{n=\frac{p}{2}}^{p-1} + \sum_{n=p+1}^{2p} + \sum_{n=2p+1}^{M} \right\} = \nonumber \\ & = & p^\alpha \{I + II + III + IV \} \nonumber
\end{eqnarray}

with the understanding that some of this sums may contain no summands (e.g. $IV = 0$ if $p \geq \frac{M}{2}$.)

Regarding $I$, if $p>3$, say, (otherwise the estimates we give are trivially true), since $\frac{n}{p-n}$ is increasing in $n$,
$$
I = \frac{1}{(p-1)^\alpha} + \frac{2^\alpha}{(p-2)^\alpha} + \dots + \frac{(\frac{p}{2}-1)^\alpha}{(\frac{p}{2}+1)^\alpha} \leq \frac{p}{2} \left\{ \frac{(\frac{p}{2})^\alpha }{ (\frac{p}{2})^\alpha }   \right\} \leq p \leq M .
$$

Also,
\begin{eqnarray}
II = \left( \frac{ p-1 }{1} \right)^\alpha +  \left( \frac{ p-2 }{2} \right)^\alpha + \dots +  \left( \frac{ \frac{p}{2} }{ \frac{p}{2}  } \right)^\alpha & \leq &  p^\alpha  \left\{ 1+ \frac{1}{2^\alpha} +  \frac{1}{3^\alpha} + \dots + \frac{1}{ \left( \frac{p}{2} \right)^\alpha}  \right\} \approx \nonumber \\
& \approx & p^\alpha \int_{1}^{\frac{p}{2}} \frac{1}{x^\alpha} dx \approx p \leq M. \nonumber
\end{eqnarray}

Regarding $III$, if $p>3$, say, (otherwise the estimates we give are trivially true),

\begin{eqnarray}
III = \left( \frac{ p+1 }{1} \right)^\alpha +  \left( \frac{ p+2 }{2} \right)^\alpha + \dots + \left( \frac{ 2p }{p} \right)^\alpha & \leq & (2p)^\alpha \left\{ 1+ \frac{1}{2^\alpha} +  \frac{1}{3^\alpha} + \dots + \frac{1}{ p^\alpha}  \right\} \lesssim \nonumber \\
& \lesssim &  p^\alpha \int_{1}^{p} \frac{1}{x^\alpha} dx \approx p \leq M. \nonumber
\end{eqnarray}

And finally for $IV$, since $\frac{n}{n-p}$ is a decreasing function of $n$, assuming $2p < M$ (otherwise $IV =0$),

$$
IV = \left( \frac{ 2p+1 }{p+1} \right)^\alpha +  \dots + \left( \frac{ M }{M-p} \right)^\alpha \leq (M-2p) \left( \frac{2p}{p} \right)^\alpha \lesssim M
$$

Now note that there are $M$ possible choices for points $a$ with first coordinate $\frac{1}{p}$, so, summing over them, and taking into account that $\displaystyle{ \sum_{p=1}^{M} p^\alpha \approx \int_{1}^{M} x^\alpha dx \approx M^{1+\alpha}}$, and doing the same reasoning for the interactions of $a$ with its column, we finally get

$$
\sum_{\substack{ a,a' \in A_N \\ a \not=a' } } {|a-a'|}^{-\alpha} \lesssim M^{4+\alpha}.
$$

With respect to the number of Euclidean distances determined by the family of sets $A_N$ and its Hausdorff dimension, let us fix $\varepsilon >0$. Consider the set $ D_{M,\varepsilon} = \left\{ \frac{1}{n} : n = M-M^{1-\frac{\varepsilon}{4}} +1, \dots , M \right\}$. Notice that $M-M^{1-\frac{\varepsilon}{4}} > \frac{M}{2}$ for sufficiently large $M$, so that the distances between any two consecutive points in $D_{M,\varepsilon}$ are all comparable with absolute constants to $\frac{1}{M^2}$. Hence, the set $ C_{N,\varepsilon} = D_{M,\varepsilon} \times D_{M,\varepsilon} \subset A_N$ has cardinality $N^{1-\frac{\varepsilon}{2}}$, since $M^2=N$, and is a Delone set. Consequently, $\dim_{\mathcal{H}}(A_N)=2$ and, if we assume the Falconer distance conjecture, by Theorem \ref{iosevichlaba} we get that

$$
\# \Delta(A_N) \geq  \# \Delta(C_{N,\varepsilon}) \gtrapprox C_\varepsilon \ N^{1-\frac{\varepsilon}{2}} \geq C'_\varepsilon \ N^{1-\varepsilon}.
$$

\end{proof}






\begin{remark}
When we define Hausdorff $\alpha$-adaptability and Minkowski $\alpha$-adaptability in the discrete setting, it is clear that some sets will have ``lower dimension" than they should for a ``stupid" reason. Namely, if we pick e.g. a $1$-separated Delone set $A_N \subset [0,N^{\frac{1}{d}}]^d$ with $\# (A_N) = N$ and add to it a few points very far away (which are also $1$-separated among themselves), calling the resulting set $S_N$, then the cardinality has essentially not changed at all, but the diameter has increased enormously, so that \eqref{diameterlim} is no longer satisfied with $\alpha = d$, but is only satisfied for much smaller values of $\alpha$. Similarly, for Hausdorff $\alpha$-adaptability, the interaction of the added points $a' \in S_N \setminus A_N$ among themselves and with the points in $A_N$ is very small, but again the diameter has increased enormously, so \eqref{OffDiagonalTermInEnergyIntegral} would no longer be satisfied with $\alpha = d$, but would only be satisfied for much smaller values of $\alpha$.

Since our aim is to apply all this machinery to the Erd\H os distance conjecture, where we can always substitute a set of cardinality $N$ by subsets of cardinality $N^{1-\varepsilon}$, for all $\varepsilon >0$ sufficiently small, it is only natural that we should allow for such small outliers (meaning $S_N \setminus A_N$) to be removed from the set. However, intuition here is likely to be misleading, since, for large $N$, $N^{1-\varepsilon}$ is much smaller than any constant fraction of $M$ (i.e. fractions of the type $\frac{M}{1000}$), so we are allowing to throw out ``most" of the set. So what seemingly is the behaviour of ``most" of the set, suddenly is completely irrelevant. The example from Theorem \ref{ExampleMinkowskiDimensionDifferentToHausdorffDimension} highlights this point, in what we believe to be a counter-intuitive instance.

A consequence of the example from Theorem \ref{ExampleMinkowskiDimensionDifferentToHausdorffDimension} is that a family of sets which is not Hausdorff $\alpha$-adaptable in $\mathbb{R}^d$ for any $\alpha>0$, can contain a family of subsets which is Hausdorff $\alpha$-adaptable for much larger $\alpha$, even $\alpha=d$, i.e. ``full" dimension! Admittedly, this is most disturbing from the viewpoint of a ``robust" theory of dimension per se and is not at all analogous to the continuous case. In order to fix this ``inconsistency" we needed to allow for ``large subsets" in the definition of discrete Hausdorff dimension. However, this is indeed an advantage for the applications of the machinery to the Erd\H os distance conjecture (which is a main point of the machinery), as we have seen in the example from Theorem \ref{ExampleMinkowskiDimensionDifferentToHausdorffDimension}, since we may verify the Erd\H os distance conjecture for a family of sets via such a ``most disturbing" family of subsets.

\end{remark}

\vskip.125in

We will now construct a family of $1$-separated finite sets $A_N \subset \mathbb{R}^d$, with $\#(A_N)=N$ so that they are not Hausdorff $\alpha$-adaptable for any $\alpha \geq 1$ in the plane. However, we do not want the family $A_N$ to be not Hausdorff $\alpha$-adaptable for any $\alpha \geq 1$ for the ``simple" aforementioned reason that most of the set is Hausdorff $\alpha$-adaptable for some $\alpha \geq 1$, but there is a small cluster (or even a single point) located very far away from the rest of the set which makes the diameter of the set huge without essentially increasing the cardinality of the main cluster of the set. Since for the Erd\H os distance problem we are allowed to remove from a set of cardinality $N$ subsets of cardinality $N-N^{1-\varepsilon}$, for $\varepsilon>0$ arbitrarily small, the example should be such that no subsets $B_N$ of these $A_N$ with $\# (B_N) \approx N^{1-\varepsilon}$, for $\varepsilon>0$ very small, are Hausdorff $\alpha$-adaptable for any $\alpha \geq 1$. In other words, we want that $\dim_{\mathcal{H}}(A_N) \leq 1$.

\begin{theorem}\label{CantorSetSmallDimension}
There exists a family a family of $1$-separated finite sets $A_N \subset \mathbb{R}^d$, with cardinality of $A_N = \#\{A_N \} =N $, so that $${({[diam(A_N)]}^{-1}A_N)}_{\frac{1}{4diam(A_N)}} \subset {[0,1]}^d$$ is a family of nested sets, but the family $A_N$ is not Hausdorff $\alpha$-adaptable for any $\alpha \geq \frac{d}{2}$. Moreover, given any $\varepsilon >0$ sufficiently small, if we consider any family $B_N \subset A_N$ with $\# (B_N) \geq C_{\varepsilon} N^{1-\varepsilon}$, then the family $B_N$ is also not Hausdorff $\alpha$-adaptable for any $\alpha \geq \frac{d}{2}$. In other words, $\dim_{\mathcal{H}}(A_N) \leq \frac{d}{2}$.
\end{theorem}

\begin{proof}
The philosophy is to mimic the construction of a Cantor set $C$ of small Hausdorff dimension $d_0$, and observe that any subset of $C$ has Hausdorff dimension $\leq d_0$. However, while this philosophy (of subsets having smaller Hausdorff dimension than the original set) works for the example we are about to construct (due to self-similarity), we already saw that it fails completely in the general case (see Theorem \ref{ExampleMinkowskiDimensionDifferentToHausdorffDimension}.) For simplicity we perform the construction in the plane.

For the construction of the Cantor set, we follow the notation and setup in \cite{mattila}. Let $0< \lambda < \frac{1}{2}$. Denote $I_{0,1} = [0,1]$, and let $I_{1,1}$ and $I_{1,2}$ be the intervals $[0,\lambda]$ and $[1-\lambda, 1]$ respectively. For each already given interval, continue the process of selecting two subintervals. If the intervals $I_{k-1,1}, \dots , I_{k-1,2^{k-1}}$ have already been defined, then define $I_{k,1}, \dots , I_{k,2^k}$ by deleting from the middle of each $I_{k-1,j}$ an interval of length $(1-2\lambda) \; diam(I_{k-1,j}) = (1-2\lambda) \lambda^{k-1}$. Thus, $length(I_{k,j}) = \lambda^{k}$.

Then define $\displaystyle  C_1 (\lambda) = \bigcap_{k=0}^{\infty} \bigcup_{j=1}^{2^k} I_{k,j} $, and $\displaystyle  C (\lambda) =  C_1 (\lambda) \times  C_1 (\lambda)$. Then $C (\lambda)$ satisfies the open set condition and $dim_{\mathcal{H}} ( C(\lambda) ) = \frac{\log (4)}{\log( \frac{1}{\lambda}) }$, which suggests that we should look for $\lambda < \frac{1}{4}$.

Consider now the previous construction up to step (or generation) $M$, for large $M$, i.e. $k=M$. Place a point in the center of each of the $N=4^M$ squares (or at any other distinguished point of the squares, but the same distinguished point for all squares, i.e. the center, the upper left corner, etc.), and set that to be $\widetilde{A_N}$. Then the minimum distance among two points in $\widetilde{A_N}$ is $\delta = (1-\lambda) \lambda^{M-1}$. Hence, in order to make the set $1$-separated, we define $A_N = \frac{1}{\delta} \widetilde{A_N}$. Consequently, $diam(A_N) \approx \frac{1}{\lambda^{M-1}}$.

Then $${({[diam(A_N)]}^{-1}A_N)}_{\frac{1}{4diam(A_N)}} \subset {[0,1]}^d$$ is a family of nested sets as long as $\lambda$ is sufficiently small (elementary calculations yield that $\lambda \lessapprox 0.1329 \dots$ is enough, although if we had considered $\frac{1}{2diam(A_N)}$ instead of $\frac{1}{4diam(A_N)}$ a larger $\lambda$ would also have worked.)

Then \eqref{diameterlim} is satisfied by $A_N$ if and only if $diam (A_N) \approx \frac{1}{\lambda^{M-1}} \lesssim 4^{\frac{M}{\alpha}}$, which is turn is true iff $\left( 4^{\frac{1}{\alpha}} \lambda \right)^M  \frac{1}{\lambda} \gtrsim 1$, which is false for $\alpha \geq 1$, since for such $\alpha$, $\left( 4^{\frac{1}{\alpha}} \lambda \right)^M \to 0 $ as $M$ (and hence $N$) $\to \infty$ (recall that $\lambda < \frac{1}{4}$.)

Now given $\varepsilon >0$ very small, consider a corresponding family $B_N \subset A_N$ with $\#(B_N) \geq C_{\varepsilon} N^{1-\varepsilon}$. Fix $\alpha \geq 1$. Since $A_N$ does not satisfy \eqref{diameterlim}, we have that $diam (A_N) > > N^{\frac{1}{\alpha}}$. In order to have any chance of $B_N$ satisfying \eqref{diameterlim}, the diameter of $B_N$ should be much smaller than that of $A_N$. Let us think in terms of starting with $A_N$ and removing successively points in order to get to $B_N$. There are only 2 procedures to reduce the diameter of $A_N$ in a substantial way by removing points from $A_N$.

The first such procedure (let us call it P1) to reduce the diameter of $A_N$ in a substantial way by removing points from $A_N$ is to at least remove 3 of the 4 squares of the form $I_{1,j} \times I_{1,k}$ and all their children. Let us call the operation of removing the 3 siblings of a given square of sidelength $2^{-k}$ (and all their descendants), an operation P (for pruning) at scale $k$. In that manner (i.e. after an operation P at scale $k=1$), the diameter of $A_N$ gets reduced by a factor of $\lambda$, and the number of points changes from $N$ to $\frac{N}{4}$. (Otherwise, if any two points contained in two different squares of the form $I_{1,j} \times I_{1,k}$ survive, the diameter of the subset of $A_N$ thus chosen is comparable to that of $A_N$.)

So, if there is any hope of $B_N$ satisfying \eqref{diameterlim}, then $B_N$ should be obtained from $A_N$ by performing an operation P at scale $k=1$, and then performing another operation P at scale $k=2$ on the surviving squares, and so on until a generation $k=L$, and then possibly removing some more points, (but not an operation of type P at generation $L+1$.) Since on the right hand side of \eqref{diameterlim} we have the number of points of the set in question, and unless we remove 3 squares (and their children) out of 4 from a given generation (i.e. we perform an operation of type P), the diameter does not decrease substantially, the best possible case given that we already performed operations P at scales $1$ through $L$ and we are not performing any further operations P, is not to remove any further points at all from the surviving squares after those consecutive $L$ operations P, in order to maximize the right hand side, once the diameter of $B_N$ is essentially fixed after those $L$ operations. This reasoning describes the candidate for $B_N$ with best chances of satisfying \eqref{diameterlim}, let us call it $\widetilde{B_N}$, in the sense that if any $B_N$ with the required conditions satisfies \eqref{diameterlim}, then so does $\widetilde{B_N}$. However, $\widetilde{B_N} = A_{\widetilde{N}}$, for some large $\widetilde{N}$ (that can be calculated explicitly, since $\widetilde{N} = \#(B_N)$), so $\widetilde{B_N}$ does not satisfy \eqref{diameterlim}, by the reasoning done for the sets $A_N$.

The reader may care to check that, indeed, for any $\alpha \geq 1$, the bound for ${\cal I}_{\alpha}(A_N)$ in equation \eqref{OffDiagonalTermInEnergyIntegral} is not satisfied, nor is it satisfied for any $B_N$ as in the statement of the Theorem.

There is however, a second procedure (let us call it P2) to reduce the diameter of $A_N$ in a substantial way by removing points from $A_N$. Namely, leaving the diameter of $A_N$ as it is, but increasing the minimum separation of the points, so that the resulting set, when rescaled to be $1$-separated, has smaller diameter.

The reader may rightfully point out that indeed these two procedures (P1 and P2) could be combined. We will deal with that possibility momentarily. Let us focus for the time being on P2. If we leave the diameter of $A_N$ untouched, but we want to increase the minimum separation between points in a substantial way, the only way to do that is to prune at the smallest scale and then move upwards in the scales. I.e. for each group of sibling squares at scale $k$, remove 3 of the 4 siblings. Let us call this operation an operation $P'$ at scale $k$. After such an operation $P'$ at scale $M$, the minimum separation between points in $A_N$ gets increased by a factor of $\frac{1}{\lambda}$, and the number of points changes from $N$ to $\frac{N}{4}$. As with P1, by a similar reasoning, the candidates for $B_N$ with best chances of satisfying \eqref{diameterlim} (let us call any of them $(B_N)'$) are the result of performing consecutively $L$ operations $P'$ and not removing any further point from $A_N$. Notice now that, after rescaling, except for the fact that the points chosen in any of the squares are not the center of the squares (or the same distinguished point in each of the squares), any such $(B_N)' = A_{N'}$, for some large $N'$ (again with $N' = \#(B_N)$), actually, $N' = 4^{M - L}$.

However it is immaterial where we place the actual points of a given set $(B_N)'$ inside each square of generation $M-L$ in the Cantor set, provided we place one point per square of generation $M-L$. To be sure, let us denote any two squares of generation $M-L$ in the Cantor set by $Q$ and $Q'$. Then for any pair of points $x,y \in Q$ and any pair of points $x',y' \in Q'$, we have that $|x-y| \approx |x'-y'|$, with comparability constants that only depend on $\lambda$ and not on $Q$ or $Q'$. Hence, if any statement regarding Hausdorff or Minkowski $\alpha$-adaptability (or dimension) of the type $\geq, \leq ,= $ (something) is true for any particular $(B_N)'$, it is simultaneously true for all such $(B_N)'$ and for $A_{N'}$. So the reasoning for P2 gets reduced to the reasoning for P1.

In a similar fashion, combining procedures P1 and P2 would yield (up to allocation of points inside each square of the smallest surviving generation) another rescaled version of $A_N$ and the same conclusion applies.

\end{proof}





The example from Theorem \ref{CantorSetSmallDimension} can be worsened to ``Hausdorff dimension $0$" as our next Theorem shows.

\begin{theorem}\label{CantorTypeSetDimensionZero}
There exists a family of sets $A_N \subseteq [0,1]^d$, with $\#(A_N)=N$, so that $(A_N)_{\delta_N}$ is a nested family of sets for some $\delta_N >0$, but so that it is not Hausdorff $\alpha$-adaptable for any $\alpha>0$. Moreover, for any $\alpha >0$, and for any family of subsets $B_N \subseteq A_N$ with $\#(B_N) \geq C_{\varepsilon} N^{1-\varepsilon}$, for sufficiently small $\varepsilon$, $B_N$ is not Hausdorff $\alpha$-adaptable. In other words, $\dim_{\mathcal{H}} (A_N)=0$.
\end{theorem}

\begin{proof}
The idea is to build a Cantor type set with decreasing proportions of ``surviving intervals" as the number of generation increases. For simplicity we perform the construction in the plane. The construction and the proof is very similar to that of Theorem \ref{CantorSetSmallDimension}.

We somewhat follow the notation and setup in \cite{mattila}. Let $0< \lambda < \frac{1}{4}$. Denote $I_{0,1} = [0,1]$, and let $I_{1,1}$ and $I_{1,2}$ be the intervals $[0,\lambda]$ and $[1-\lambda, 1]$ respectively. For each already given interval, continue the process of selecting two subintervals. If the intervals $I_{k-1,1}, \dots , I_{k-1,2^{k-1}}$ have already been defined, then define $I_{k,1}, \dots , I_{k,2^k}$ by keeping from each $I_{k-1,j}$ two intervals of length $f_k := \frac{\lambda}{2^{k-1}}$ times the length of $I_{k-1,j}$ with the same endpoints as $I_{k-1,j}$ (the notation $f_k$ stands for ``factor at scale $k$".) Thus, $length(I_{k,j}) = \frac{\lambda^{k}}{2^{\frac{k(k-1)}{2}}}$. Notice that $f_k$ decreases as $k$ increases.

Then define $\displaystyle  C_1 (\lambda) = \bigcap_{k=0}^{\infty} \bigcup_{j=1}^{2^k} I_{k,j} $, and $\displaystyle  C (\lambda) =  C_1 (\lambda) \times  C_1 (\lambda)$.

Since at stage $M$ of the previous construction there are $N=4^M$ squares of sidelength $\frac{\lambda^{M}}{2^{\frac{M(M-1)}{2}}}$, an easy calculation yields that $\dim_{\mathcal{H}} C (\lambda) =0$. Let us take a point in each of the aforementioned $N=4^M$ squares and let the resulting set be $\widetilde{A_N}$.

Let us briefly remark that
it is immediate from the continuous case calculations that $\widetilde{A_N}$ is not Hausdorff $\alpha$-adaptable for any $\alpha >0$. Namely, fix $\alpha >0$ and take $C_N := \left(  \widetilde{A_N}  \right)_{N^{- \; \frac{1}{\alpha}}}$. Then, as in the proof of Theorem \ref{HausdorffAlphaAdaptableSatisfyErdosIfAlphaAccordingToFalconerVersion2WhenSetsAreNested}, $C_N \to C(\lambda)$ in the Hausdorff metric, and then if $\widetilde{A_N}$ were Hausdorff $\alpha$-adaptable, the energy integral $I_\alpha (\mu_{A_N}) \leq C < \infty$ for all $N$. By taking a subsequence, we could assume that $\mu_{A_N} \rightharpoonup \mu_0$, in the sense of weak-$\ast$ convergence, and then $supp (\mu_0) \subseteq C(\lambda)$. Then $I_\alpha (\mu_0) \leq C$, so that $\dim_{\mathcal{H}} \left(  C(\lambda)  \right) \geq \dim_{\mathcal{H}} \left(  supp (\mu_0) \right) \geq \alpha$, which would be a contradiction. However we prefer to do direct calculations in order to show that \eqref{OffDiagonalTermInEnergyIntegralWithLessapprox} is also not satisfied.

The minimum separation between points in the set $\widetilde{A_N}$ is $\displaystyle{\approx \left( 1- \frac{4\lambda}{2^M}   \right) \frac{\lambda^{M-1}}{2^{ \frac{(M-1)(M-2)}{2}  }} }$, so in order to make the set $\widetilde{A_N}$ 1-separated, we have to rescale by the inverse of the minimum separation between points which is $$\approx \frac{2^{ \frac{(M-1)(M-2)}{2} } }{\lambda^{M-1}} = diam (A_N), $$ denoting by $A_N$ such a rescaling of $\widetilde{A_N}$.

If the family of sets $A_N$ were Hausdorff $\alpha$-adaptable, for some $\alpha>0$, then we would need that $diam (A_N) \lesssim N^{\frac{1}{\alpha}} = 4^{\frac{M}{\alpha}}$, by \eqref{diameterlim}. But this is equivalent to

$$
2^{ \frac{(M-1)(M-2)}{2} }  \leq C 4^{\frac{M}{\alpha}} \lambda^{M-1},
$$

which in turn, taking logarithms, is equivalent to

$$
\frac{(M-1)(M-2)}{2} \leq \frac{2M}{\alpha} + C_1 M + C_2
$$

for some constants $C_1, C_2$, which is impossible if $M \to \infty$, for any $\alpha >0$.

Now fix $\varepsilon >0$ sufficiently small and assume we have a sequence of subsets $B_N \subset A_N$ with $\#(B_N) \geq C_{\varepsilon} N^{1-\varepsilon}$. Let us fix some $\alpha>0$. If the family $B_N$ has any chance of being Hausdorff $\alpha$-adaptable, then the diameter of $B_N$ should be considerably smaller (after rescaling $B_N$ to be $1$-separated) than that of $A_N$, since by the proof of $A_N$ not being Hausdorff $\alpha$-adaptable, we know that $diam(A_N) > > N^{\frac{1}{\alpha}}$.
Let us again think in terms of removing points from $A_N$ in order to get to $B_N$.
As in Theorem \ref{CantorSetSmallDimension}, there are only 2 procedures to substantially reduce the diameter of the resulting set starting from $A_N$.

The first procedure (P1), consists again of removing 3 of the 4 squares of the form $I_{1,j} \times I_{1,k}$ and all their children (i.e. performing an operation P at scale $k=1$), and then repeating the same operation with 3 of the 4 surviving squares of generation 2, and so on, repeating the operation P exactly for the first $L$ scales. Once this operation has been performed exactly $L$ times, the diameters of the possible subsets $B_N$ (i.e. if no further operation P is performed) are all comparable, and hence the $B_N$ with best possible chances is the one with most points, i.e. the set with no further points removed after those $L$ operations P. Since each operation P divides the number of points by 4, we have that $\frac{N}{4^L} = \#(B_N) \geq C_{\varepsilon} N^{1-\varepsilon}$.

The second procedure (P2), consists again of removing of removing 3 of the 4 siblings for each group of sibling squares at scale $k$ (let us again call this operation an operation $P'$ at scale $k$), starting from the smallest scale and moving up in the scales. Each operation $P'$ divides the number of points by 4, as with operation P. However, since in our present case the factors $f_k$ are not constant (as they were in Theorem \ref{CantorSetSmallDimension}), but they are decreasing in $k$, now the operation $P'$ is substantially more efficient than the operation $P$ in terms of reducing the diameter of the set in question (after rescaling the set so that it is $1$-separated.)

Consequently, the candidate for $B_N$ with best chances of being Hausdorff $\alpha$-adaptable (let us call it $(B_N)'$) is the result of performing the procedure P2 from the smallest scale, moving up the scales, exactly $L$ times and not removing any further point from $A_N$. But, after rescaling so that $(B_N)'$ becomes $1$-separated, as in Theorem \ref{CantorSetSmallDimension}, $(B_N)' = A_{N'}$ for a certain large $N'$ ($N' = 4^{M-L}$), except for the location of the points inside each of the squares of the smallest scale (those of generation $M-L$). As in Theorem \ref{CantorSetSmallDimension}, the location of the points inside each of the squares of generation $M-L$ is immaterial for Minkowski or Hausdorff $\alpha$-adaptability (or dimension) purposes, so we can assume without loss of generality that $(B_N)'$ is really $= A_{N'}$, which we already know is not Hausdorff $\alpha$-adaptable. So we get that $\dim_{\mathcal{H}}(A_N)=0$.









\end{proof}

As a concluding remark, notice that this paper highlights, among other things, that the notion of Hausdorff dimension (even in the continuous case) contains much more information than just the size of the sets, since, after all, all the families of sets we described have the same size (namely $N$.) Hausdorff dimension is more about ``electrostatics" (how different charges are positioned relatively to one another) than about size. (The case of $\mathbb{R}^3$ and $\alpha=1$ is indeed classical electrostatics and the energy integral we considered is the energy of the system of charges.)

\bibliographystyle{alpha}

\end{document}